\documentclass[11pt]{article}

\usepackage[utf8]{inputenc}
\usepackage[T1]{fontenc}
\usepackage[a4paper, margin=2.7cm]{geometry}
\usepackage{amsmath, amsthm, amsfonts, amssymb}
\usepackage{mathrsfs}           

\usepackage[colorlinks=true,linkcolor=blue,citecolor=blue,urlcolor=blue,breaklinks]{hyperref}

\usepackage{graphicx}
\usepackage{amsmath}
\usepackage{latexsym}
\usepackage{amssymb}
\usepackage{amsbsy}
\usepackage{dsfont}
\usepackage{cite}
\usepackage{subfigure}
\def\be{\begin{equation}}
\def\ee{\end{equation}}
\def\bea{\begin{eqnarray}}
\def\beas{\begin{eqnarray*}}
\def\ba{\begin{array}{l}\displaystyle}
\def\eea{\end{eqnarray}}
\def\eeas{\end{eqnarray*}}
\def\ea{\end{array}}
\def\R{\mathbb{R}}

\def\nuu{s}

\newtheorem{thm}{Theorem}[section]

\theoremstyle{definition}

\theoremstyle{remark}
\newtheorem{remark}[thm]{Remark}
\theoremstyle{example}

\title{High order semi-implicit multistep methods for time dependent partial differential equations
}

\author{Giacomo Albi\thanks{Computer Science Department, University of Verona, Italy}        \and
        Lorenzo Pareschi\thanks{Mathematics and Computer Science Department, University of Ferrara, Italy} 
}

\date{December 31, 2019}
\AtEndDocument{\bigskip{\footnotesize
		\noindent
		\textsc{G. Albi, str. Le Grazie 15, 37134, Verona, Italy.}
		\textit{E-mail address}: \texttt{\href{ giacomo.albi@univr.it}{giacomo.albi@univr.it}}
		\par
		\addvspace{\medskipamount}
		\noindent
		\textsc{L. Pareschi, via Machiavelli 30, 44121, Ferrara, Italy.}
		\textit{E-mail address}: \texttt{\href{lorenzo.pareschi@unife.it}{lorenzo.pareschi@unife.it}}
}}

\begin{document}
	
\maketitle

\begin{abstract}
We consider the construction of semi-implicit linear multistep methods which can be applied to time dependent PDEs where the separation of scales in additive form, typically used in implicit-explicit (IMEX) methods, is not possible. As shown in \cite{BFR} for Runge-Kutta methods, these semi-implicit techniques
give a great flexibility, and allows, in many cases, the construction of
simple linearly implicit schemes with no need of iterative solvers. In this work we develop a general setting for the construction of high order semi-implicit linear multistep methods and analyze their stability properties for a prototype linear advection-diffusion equation and in the setting of strong stability preserving (SSP) methods.
Our findings are demonstrated on several examples, including nonlinear reaction-diffusion and convection-diffusion problems.
\end{abstract}
\paragraph{keywords:}{semi-implicit methods \and implicit-explicit methods \and multistep methods \and strong stability preserving \and high order accuracy}


\section{Introduction}
\label{intro}
In many situations we have to deal with dynamical system arising from the spatial discretization (or more in general discretization in the phase space) of time dependent partial differential equations. For problems where the various terms have different time scales which can be easily separated the resulting system of ordinary differential equations has the form
\be
\label{eq:PAR}
\frac{d u(t)}{d t} = f(t,u(t))+\frac1{\varepsilon}g(t,u(t)),
\ee
with $\varepsilon>0$ a small parameter emphasizing the stiffness in the system. In \eqref{eq:PAR} the solution $u(t)$ is a vector in $\R^m$ which initially satisfies the condition $u(0)=u_0$. Typically, the term $f$ contains some nonlinearity or complexity that we do not want to integrate implicitly, whereas the term $g/\varepsilon$ is stiff and requires an implicit integration. For systems of the type \eqref{eq:PAR} implicit-explicit (IMEX) schemes are nowadays a very popular choice \cite{ARS, Ascher2, CK, PR}.  

However, in some cases, this separation is not possible and more in general we are lead to a dynamical system in the form 
\be
\label{eq:SYS}
\frac{d u(t)}{d t} = {\mathcal H}\left(t,u(t),\frac{u(t)}{\varepsilon}\right)
\ee 
where the right hand side has a stiff dependence only on the last argument.

Similarly to \eqref{eq:PAR} for problem \eqref{eq:SYS} it is highly desirable to construct a numerical method based on evaluating implicitly only the stiff component $u(t)/\varepsilon$ by keeping explicit the non stiff one in order to reduce the computational complexity of a fully implicit solver. Following \cite{BFR} we shall call them semi--implicit methods, since they can be used in a more general context than implicit-explicit methods.

A simple first order method which realizes the above idea reads as follows
\be
\label{eq:SYS1}
u^{n+1} = u^n + \Delta t {\mathcal H}\left(t^{n+1},u^{n},\frac{u^{n+1}}{\varepsilon}\right).
\ee 
There are many circumstances in which this semi-implicit approach leads to considerable advantages, for example if the function  ${\mathcal H}$ is linear with respect to the stiff component the numerical solution can be computed solving only linear systems of equations\cite{BBMRV}. 

The extension of this simple idea to higher order, however, is not straightforward. 
Here, following the approach recently introduced in \cite{BFR} for Runge-Kutta methods we construct high order semi-implicit discretization based on the use of linear multistep methods. 
To this aim, by setting $v(t)=u(t)/\varepsilon$ in \eqref{eq:SYS} we obtain the equivalent formulation as a partitioned system
\bea
\nonumber
\frac{d u(t)}{d t} &=& {\mathcal H}\left(t,u(t),v(t)\right)\\[-.2cm]
\label{eq:GPS}
\\[-.2cm]
\nonumber
\varepsilon\frac{d v(t)}{d t} &=& {\mathcal H}\left(t,u(t),v(t)\right),
\eea
where $v(0)=u_0/\varepsilon$.
Thus, the formal equivalence among \eqref{eq:SYS} and \eqref{eq:GPS} allows us to adopt IMEX techniques for partitioned systems to more general cases \cite{BFR,BP}. We refer to \cite{BFR} for a detailed discussion on the equivalence between the various forms of system usually treated with IMEX methods and to \cite{ARS, CK} for general references on IMEX methods based on Runge-Kutta schemes. A large literature in this direction has been devoted to the construction of IMEX Runge-Kutta schemes satisfying the asymptotic-preserving (AP) property in the case of hyperbolic problems \cite{BP, BPR17, Jin, PR} and for kinetic equations \cite{BPR,JF,DP,DP15}. 
For the case of IMEX linear multistep methods, we refer to \cite{Ak1, Ak2, AHP1, Ascher2, FHV, HR, HRS, RSSZ1, RSSZ2} for results on the construction and properties of the schemes for various types of PDEs and to \cite{ADP,DPLMM} for the construction of schemes satisfying the AP property.

The rest of the manuscript is organized as follows. In the next Section we introduce the IMEX multistep methods for partitioned systems and derive the corresponding semi-implicit formulation for problem \eqref{eq:SYS}. Next in Section 3 we detail the derivation of the schemes and analyze their stability properties for a prototype linear advection-diffusion equation and in the setting of strong stability preserving (SSP) methods. Section 4 is then devoted to present several numerical applications that confirm the validity of the present approach. The manuscript ends with some conclusions in Section 5.

\section{Semi-implicit multistep methods}
In this Section, we first introduce the general class of IMEX linear multistep schemes for partitioned systems together with some preliminary definitions. Next we recall some general results on the order conditions and, subsequently, we apply the schemes to system in the general form \eqref{eq:SYS} to derive the corresponding semi-implicit formulation. 

\subsection{IMEX linear multistep methods for partitioned systems}

Let us consider a general partitioned system in the form
\bea
\nonumber
\frac{d y(t)}{d t} &=& {\mathcal F}\left(t,y(t),z(t)\right)\\[-.2cm]
\label{eq:GPS2}
\\[-.2cm]
\nonumber
\frac{d z(t)}{d t} &=& {\mathcal G}\left(t,y(t),z(t)\right).
\eea
where $y(t)\in\R^p$ and $z(t)\in\R^q$, $p,q\geq 1$ and $y(0) = y_0$, $z(0) = z_0$.

For partitioned system in the form (\ref{eq:GPS2}) we consider schemes based on solving the first component with an explicit linear multistep method and the second with an implicit one 
\bea
\nonumber
   y^{n+1} = - \sum_{j=0}^{\nuu-1} \tilde a_j y^{n-j} + \Delta t \sum_{j=0}^{\nuu-1} \tilde b_j {\mathcal F}\left(t^{n-j},y^{n-j},z^{n-j}\right)\\[-.2cm]
   \label{eq:GPSLM}
   \\[-.2cm]
   \nonumber
   z^{n+1} = - \sum_{j=0}^{\nuu-1} a_j z^{n-j} + {\Delta t} \sum_{j=-1}^{\nuu-1} b_j  {\mathcal G}\left(t^{n-j},y^{n-j},z^{n-j}\right)
\eea
where $b_{-1} \neq 0$. Implicit methods for which $b_j=0$, $j=0,\ldots,\nuu-1$ are referred to as backward differentiation formula (BDF). 
Another important class is represented by  Adams methods, for which $\tilde a_0=-1$, $a_0=-1$, $\tilde a_j=0$, $a_j=0$, $j=1,\ldots,\nuu-1$. 

\begin{remark}
Classical systems in additive form \eqref{eq:PAR}, i.e.
\be
\label{eq:PAR2}
\frac{d u(t)}{d t} = f(t,u(t))+\frac1{\varepsilon}g(t,u(t))
\ee
with initial data $u(0)=u_0$ 
can be also written in partitioned form by defining $u(t)=y(t)+z(t)$, ${\mathcal F}(t,y(t),z(t))=f(t,u(t))$, ${\mathcal G}(t,y(t),z(t))=g(t,u(t))/\varepsilon$ and rewriting
\bea
\nonumber
\frac{d y(t)}{d t} &=& {\mathcal F}\left(t,y(t),z(t)\right)\\[-.2cm]
\label{eq:GPS3}
\\[-.2cm]
\nonumber
\frac{d z(t)}{d t} &=& {\mathcal G}\left(t,y(t),z(t)\right),
\eea 
for any initial data such that $y(0)+z(0)=u_0$.
\end{remark}

\subsection{Semi-implicit methods in predictor-corrector form}
Let us now apply the previous general formulation to the case of system \eqref{eq:SYS} in the partitioned form \eqref{eq:GPS}, i.e.
\bea
\nonumber
\frac{d u(t)}{d t} &=& {\mathcal H}\left(t,u(t),v(t)\right)\\[-.2cm]
\label{eq:GPSbis}
\\[-.2cm]
\nonumber
\frac{d v(t)}{d t} &=& {\mathcal H}\left(t,u(t),v(t)\right),
\eea
where in order to simplify notations, we removed the dependence of the parameter $\varepsilon$ in the second
argument, keeping in mind that this dependence is stiff.

We obtain the semi-implicit multistep solver
\begin{equation}\begin{aligned}\label{eq:semiGPS}
   u^{n+1} = - \sum_{j=0}^{\nuu-1} \tilde a_j u^{n-j} + \Delta t \sum_{j=0}^{\nuu-1} \tilde b_j {\mathcal H}\left(t^{n-j},u^{n-j},v^{n-j}\right)\cr
   v^{n+1} = - \sum_{j=0}^{\nuu-1} a_j v^{n-j} +{\Delta t} \sum_{j=-1}^{\nuu-1} b_j  {\mathcal H}\left(t^{n-j},u^{n-j},v^{n-j}\right)
\end{aligned}\end{equation}
where initially we assume $v^{n-j} = u^{n-j}$, $j=0,\ldots,\nuu-1$. Let us point out that even if the scheme doubles the number of unknown the number of evaluations of the function ${\mathcal H}(t,u(t),v(t))$ is not doubled since both schemes use the same time levels. In particular, if for notation simplicity we assume the system to be autonomous and the function $\mathcal H$ to depend linearly from the second argument
\be
{\mathcal H}\left(u(t),v(t)\right)={\mathcal K}(u(t))+A(u(t))v(t),
\label{eq:lin}
\ee 
where the function $\mathcal{K}: \R^m \to \R^m$ and $A(u(t))$ is an invertible $m\times m$ matrix, 
the resulting scheme can be solved without any need of an iterative solver. In fact, the second equation in \eqref{eq:semiGPS} can be rewritten as 
\beas
 v^{n+1} =& -& \sum_{j=0}^{\nuu-1} a_j v^{n-j} +{\Delta t} \sum_{j=0}^{\nuu-1} b_j  {\mathcal H}\left(u^{n-j},v^{n-j}\right)\\
 &+&{\Delta t} b_{-1}\left({\mathcal K}\left(u^{n+1}\right)+ A(u^{n+1})v^{n+1}\right).
 \eeas
 or equivalently in explicit form
\begin{align*}
 v^{n+1} = &\left(I-{b_{-1}\Delta t} A(u^{n+1})\right)^{-1}\left(-\sum_{j=0}^{\nuu-1} a_j v^{n-j} + {\Delta t} \sum_{j=0}^{\nuu-1} b_j  {\mathcal H}\left(u^{n-j},v^{n-j}\right)\right.\cr
 &\qquad\qquad\qquad+{\Delta t} b_{-1}\left({\mathcal K}\left(u^{n+1}\right)\right)\Bigg),
\end{align*}
since $u^{n+1}$ is computed from the first equation in \eqref{eq:semiGPS}.

For semi-implicit Runge-Kutta methods (see \cite{BFR}) in the case of autonomous systems it is possible to construct the scheme in such a way that the two solutions furnished by the system at time $n+1$ coincide. At variance, for scheme \eqref{eq:semiGPS} at time level $n+1$ we will have two distinct numerical solutions $u^{n+1}$ and $v^{n+1}$ approximations of the true solution $u(t^{n+1})$ of problem \eqref{eq:SYS}. Note, however, that both these solutions are used by the scheme to advance in time. 

In order to define a unique solution, it is natural to consider the scheme \eqref{eq:semiGPS} as a predictor--corrector multistep method for the non stiff component, where the explicit scheme is used to predict $u^{n+1}$ which is then used by the implicit solver as a corrector for $v^{n+1}$. 

As an example, reverting back to the notations used initially, we can write the semi-implicit scheme for \eqref{eq:SYS} as
\bea
\nonumber
   \hat u^{n+1} &=& - \sum_{j=0}^{\nuu-1} \tilde a_j u^{n-j} + \Delta t \sum_{j=0}^{\nuu-1} \tilde b_j {\mathcal H}\left(t^{n-j},u^{n-j},\frac{u^{n-j}}{\varepsilon}\right)\\[-.2cm]
   \label{eq:GPSLMLPC}
   \\[-.2cm]
   \nonumber
   u^{n+1} &=& - \sum_{j=0}^{\nuu-1} a_j u^{n-j} + {\Delta t} \sum_{j=0}^{\nuu-1} b_j  {\mathcal H}\left(t^{n-j},u^{n-j},\frac{u^{n-j}}{\varepsilon}\right)\\
   \nonumber
   &&+\Delta t b_{-1}{\mathcal H}\left(t^{n+1},\hat u^{n+1},\frac{u^{n+1}}{\varepsilon}\right).
\eea  
The above scheme uniquely identifies the numerical solution at time $u^{n+1}$. In the following we will discuss the general order conditions for \eqref{eq:GPSLM} and present different types of semi-implicit multistep methods of various order.
\begin{remark}
As a consequence of the above predictor-corrector formulation for the non stiff component, if the implicit solver has order $p$ it is typically enough to consider an explicit solver of order $p-1$.
We emphasize that this predictor-corrector interpretation holds true only for the non stiff component, since the stiff one is treated implicitly by the resulting scheme. 
\end{remark}


\section{Order conditions, stability and derivation of the schemes}
In this Section we provide the details of the schemes that will be used in our numerical results. First we recall the order condition and then some general definitions concerning the stability properties. The stability properties for IMEX multistep methods are usually discussed in the case of additive systems for simple one-dimensional convection-diffusion problems \cite{Ascher2, FHV, HR, HRS}. We also refer to the recent stability analysis in \cite{ADP} for the case of two dimensional partitioned systems.

\subsection{Order conditions}
 For a partitioned system in the form \eqref{eq:GPSLM} an order $p$ scheme is obtained if both schemes are of order $p$, namely the following conditions are satisfied
\begin{equation}
\label{eq:GPSLMcond}
\begin{aligned}
1+\sum_{j=0}^{\nuu-1} \tilde a_j =0,&\qquad 1+\sum_{j=0}^{\nuu-1} a_j =0,\\
1-\sum_{j=1}^{\nuu-1} j \tilde a_j=\sum_{j=0}^{\nuu-1} \tilde b_j,&\qquad  1-\sum_{j=1}^{\nuu-1} j a_j  =\sum_{j=-1}^{\nuu-1} b_j,\\
\frac12+\sum_{j=1}^{\nuu-1} \frac{j^2}{2}\tilde a_j=-\sum_{j=1}^{\nuu-1} j \tilde b_j,&\qquad\frac12+\sum_{j=1}^{\nuu-1} \frac{j^2}{2}a_j=b_{-1}-\sum_{j=1}^{\nuu-1}j b_j\\
\vdots\qquad\quad\qquad&\qquad\qquad\qquad\qquad\vdots\\
\frac1{p!}+\sum_{j=1}^{\nuu-1} \frac{(-j)^p}{p!} \tilde a_j\qquad\qquad &\qquad \frac1{p!}+\sum_{j=1}^{\nuu-1} \frac{(-j)^p}{p!} a_j\\
=\sum_{j=1}^{\nuu-1}\frac{(-j)^{p-1}}{(p-1)!} \tilde b_j,&\qquad\qquad\qquad 
=\frac{b_{-1}}{(p-1)!}+\sum_{j=1}^{\nuu-1}\frac{(-j)^{p-1}}{(p-1)!} b_j.
\end{aligned}
\end{equation}
We recall that a $\nuu$-step implicit multistep scheme can achieve order $\nuu+1$, while an $\nuu$-step explicit methods has only order at most $\nuu$.
We refer to~\cite{Ascher2, HR, FHV,HW} for more details on
the order conditions for IMEX multistep schemes in the case of additive systems. Here, we remark that the order conditions for the partitioned scheme are simpler than in the case of additive schemes where there is a coupling between the explicit and the implicit solver. In addition, in the case of system \eqref{eq:GPSLMLPC} only $p-1$ order of accuracy is required by the explicit predictor solver to guarantee that an order $p$ implicit corrector step yields the desired order $p$ accuracy.

\subsection{Stability properties}
The stability properties are usually analyzed for simple one-dimensional linear problems of the form
\be
\frac{d u(t)}{d t} = i\lambda u(t) + \mu u(t),
\label{eq:add}
\ee
where $\lambda,\mu\in \R$ and $i$ is the imaginary unit. These kind of problems typically are originated by the space discretization of convection diffusion equations, hyperbolic balance laws or linear kinetic models~\cite{Ascher2, DPLMM, HR, FHV,HW}. For example, in the case of one-dimensional linear convection diffusion problems 
\[
\frac{\partial u}{\partial t} = a \frac{\partial u}{\partial x}+D \frac{\partial^2 u}{\partial x^2},
\]
discretized by standard central differences of mesh $\Delta x$ we have \eqref{eq:add} with
\[
\lambda = \frac{a}{\Delta x} \sin(k\Delta x),\qquad \mu = \frac{2D}{(\Delta x)^2}\left(\cos(k\Delta x)-1\right),
\]
where $k$ is the frequency of the corresponding Fourier mode.

In partitioned form system \eqref{eq:add} will correspond to system
\bea
\nonumber
\frac{d u(t)}{d t} &=& i\lambda u(t)+ \mu v(t)\\[-.2cm]
\label{eq:GPSstab}
\\[-.2cm]
\nonumber
\frac{d v(t)}{d t} &=& i\lambda u(t)+\mu v(t).
\eea
Note that, since the direct application of an IMEX multistep method to a system in the additive form \eqref{eq:add} is not equivalent to the application of the combination of an explicit and an implicit multistep scheme to the partitioned form \eqref{eq:GPSstab} even the resulting stability analysis will differ.

Applying the semi-implicit scheme \eqref{eq:semiGPS} to \eqref{eq:GPSstab} yields 
\bea
\nonumber
  u^{n+1} =&& - \sum_{j=0}^{\nuu-1} \tilde a_j u^{n-j} + \Delta t \sum_{j=0}^{\nuu-1} \tilde b_j \left(i\lambda u^{n-j} + {\mu} v^{n-j}\right)\\[-.2cm]
   \label{eq:GPSstabd}
   \\[-.2cm]
   \nonumber
    v^{n+1} =&& \left(\frac{1}{1-\mu b_{-1}\Delta t}\right)\left(-\sum_{j=0}^{\nuu-1} a_j v^{n-j} + {\Delta t} \sum_{j=0}^{\nuu-1} b_j \left(i\lambda u^{n-j} + \mu v^{n-j}\right)\right.\\ 
    \nonumber
 &+&{\Delta t} b_{-1} i\lambda u^{n+1}\Bigg),
\eea
where $u^{n-j}=v^{n-j}$, $j=0,\ldots,\nuu-1$.
By direct substitution of the first equation into the second we obtain the explicit form
\beas
    v^{n+1} =&& \left(\frac{1}{1-\mu b_{-1}\Delta t}\right)\left(-\sum_{j=0}^{\nuu-1} a_j v^{n-j} + {\Delta t}\sum_{j=0}^{\nuu-1} b_j \left(i\lambda v^{n-j} + \mu v^{n-j}\right)\right.\\ 
    \nonumber
 &+&\left.{\Delta t} b_{-1} i\lambda \left(-\sum_{j=0}^{\nuu-1} \tilde a_j v^{n-j} + \Delta t \sum_{j=0}^{\nuu-1} \tilde b_j \left(i\lambda v^{n-j} + \mu v^{n-j}\right)\right)\right).
\eeas
This leads to the characteristic equation
\[
\zeta^{\nuu}({1-b_{-1}z_R})+\rho(\zeta)-(z_R+z_I)\sigma(\zeta)+b_{-1}z_I\left(\tilde\rho(\zeta)-(z_R+z_I)\tilde\sigma(\zeta)\right)=0,
\]
where $z_R=\mu\Delta t$, $z_I=i\lambda\Delta t$ and
\beas
\rho(\zeta)&=&\sum_{j=0}^{\nuu-1} a_j \zeta^{s-j-1},\quad \sigma(\zeta)=\sum_{j=0}^{\nuu-1} b_j \zeta^{s-j-1},\\
\tilde \rho(\zeta)&=&\sum_{j=0}^{\nuu-1} \tilde a_j \zeta^{s-j-1},\quad 
\tilde\sigma(\zeta)=\sum_{j=0}^{\nuu-1} b_j \zeta^{s-j-1}.
\eeas
Stability corresponds to the requirement that all roots have modulus less or equal one and that all multiple roots have modulus less than one. 
 

\subsection{Derivation of the schemes}

In the sequel to simplify notations, we restrict to autonomous systems (this is always possible simply by augmenting the dimension of the system by one), i.e. the function ${\mathcal H}$ does not depend explicitly on time. Let us first point out that the simplest first order method \eqref{eq:SYS1} is common to multistep methods and Runge-Kutta methods and in the case of system \eqref{eq:GPSbis} reads
\beas
u^{n+1} &=& u^n\\
v^{n+1} &=& v^n +{\Delta t}{\mathcal H}\left({u}^{n+1},v^{n+1}\right),
\eeas
with $v^n=u^n$, which corresponds to a simple identity as explicit predictor for the non stiff component and backward Euler as implicit corrector for the stiff one.   

\paragraph{Second order methods.}
The general form of second order schemes for \eqref{eq:GPSbis} reads as follows
\beas
u^{n+1} &=& u^n + {\Delta t}{\mathcal H}\left({u}^{n},v^n\right)\\
v^{n+1} &=& \frac1{2\alpha+1}\left(4\alpha u^n -(2\alpha-1)u^{n-1}\right) +\frac{\Delta t}{2\alpha+1}\left((2\alpha+\beta){\mathcal H}\left({u}^{n+1},v^{n+1}\right)\right.\\
&&\left.+
2(1-\alpha-\beta){\mathcal H}\left({u}^{n},v^n\right)+
\beta {\mathcal H}\left({u}^{n-1},v^{n-1}\right) \right),
\eeas
with $u^n=v^n$, $u^{n-1}=v^{n-1}$ and 
the solver for the non stiff component is represented by forward Euler. Popular choices for the implicit solver are obtained for $\alpha=1/2$ and $\beta=0$ which corresponds to Crank-Nicholson, the resulting scheme will be referred to as FE-CN2, and $\alpha=1$ and $\beta=0$ corresponding to the second order BDF scheme, we will refer to this scheme as FE-BDF2. 
In the case of $\alpha=1/2$ the value of $\beta=1/8$ yields the best damping properties \cite{Ascher2}, the resulting scheme requires the additional storage of level $n-1$ and is referred to as FE-MCN2.

\begin{figure}
\begin{center}
\includegraphics[scale=0.33]{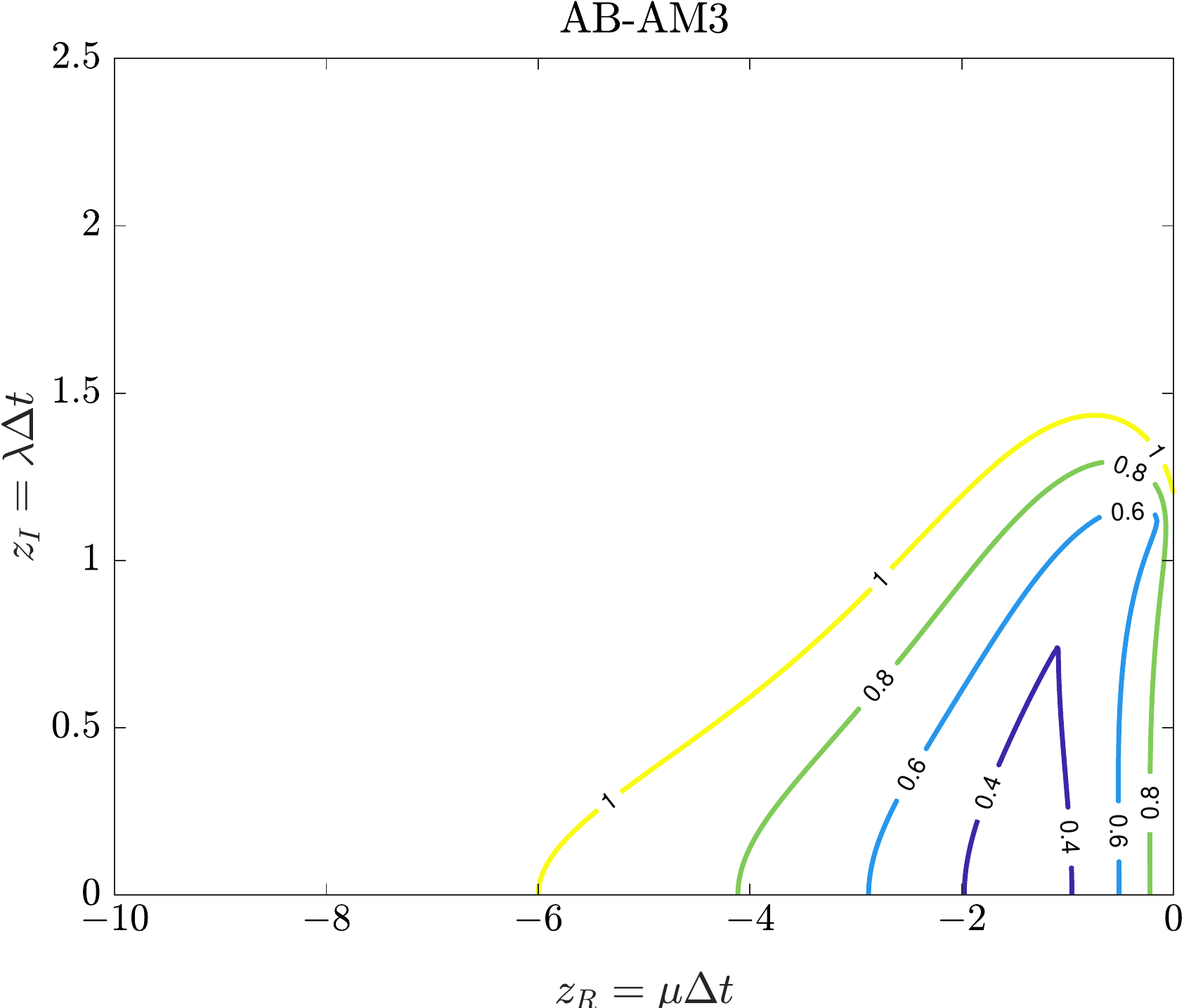}\quad
\includegraphics[scale=0.33]{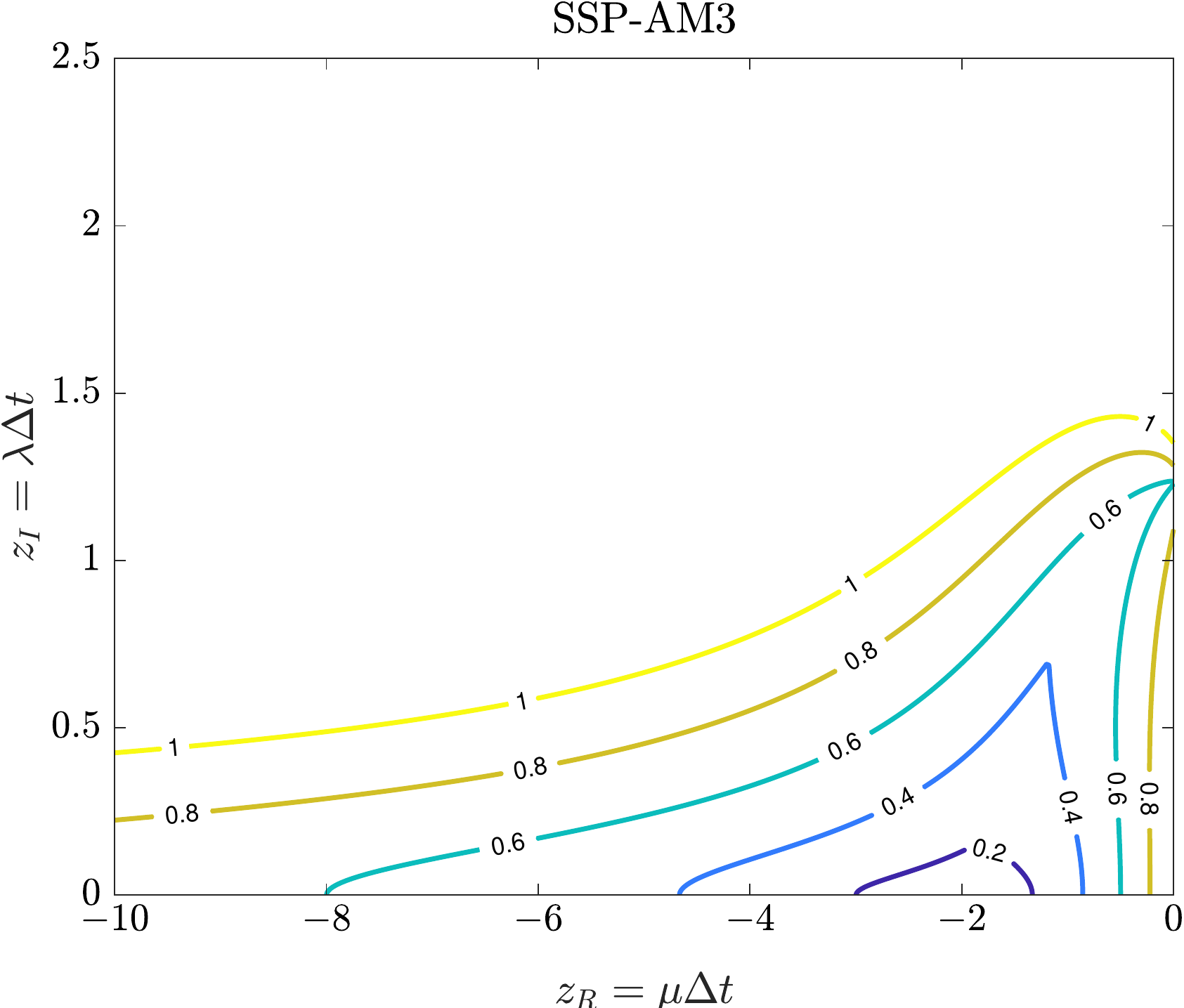}\\[+.1cm]
\includegraphics[scale=0.33]{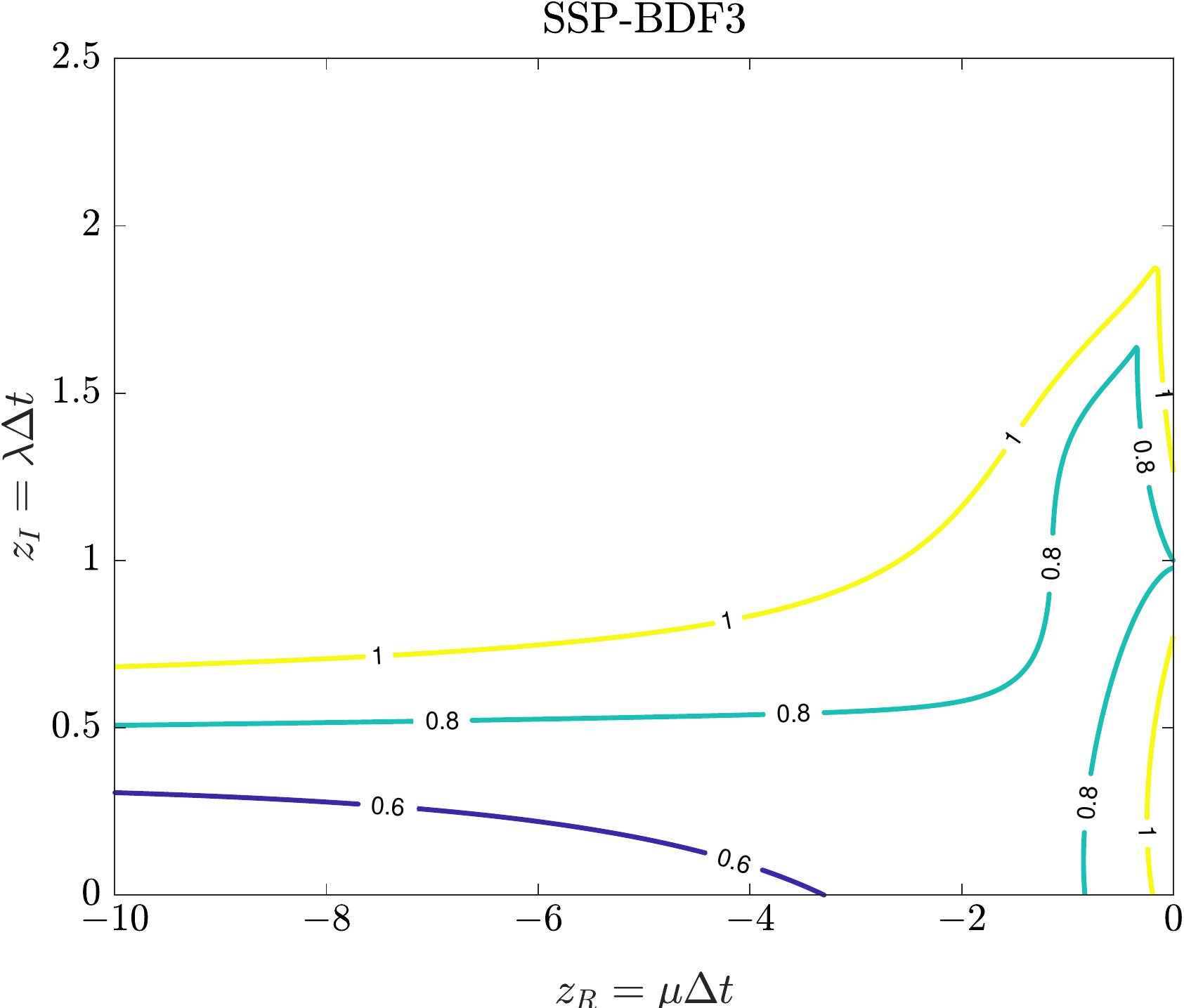}\quad
\includegraphics[scale=0.33]{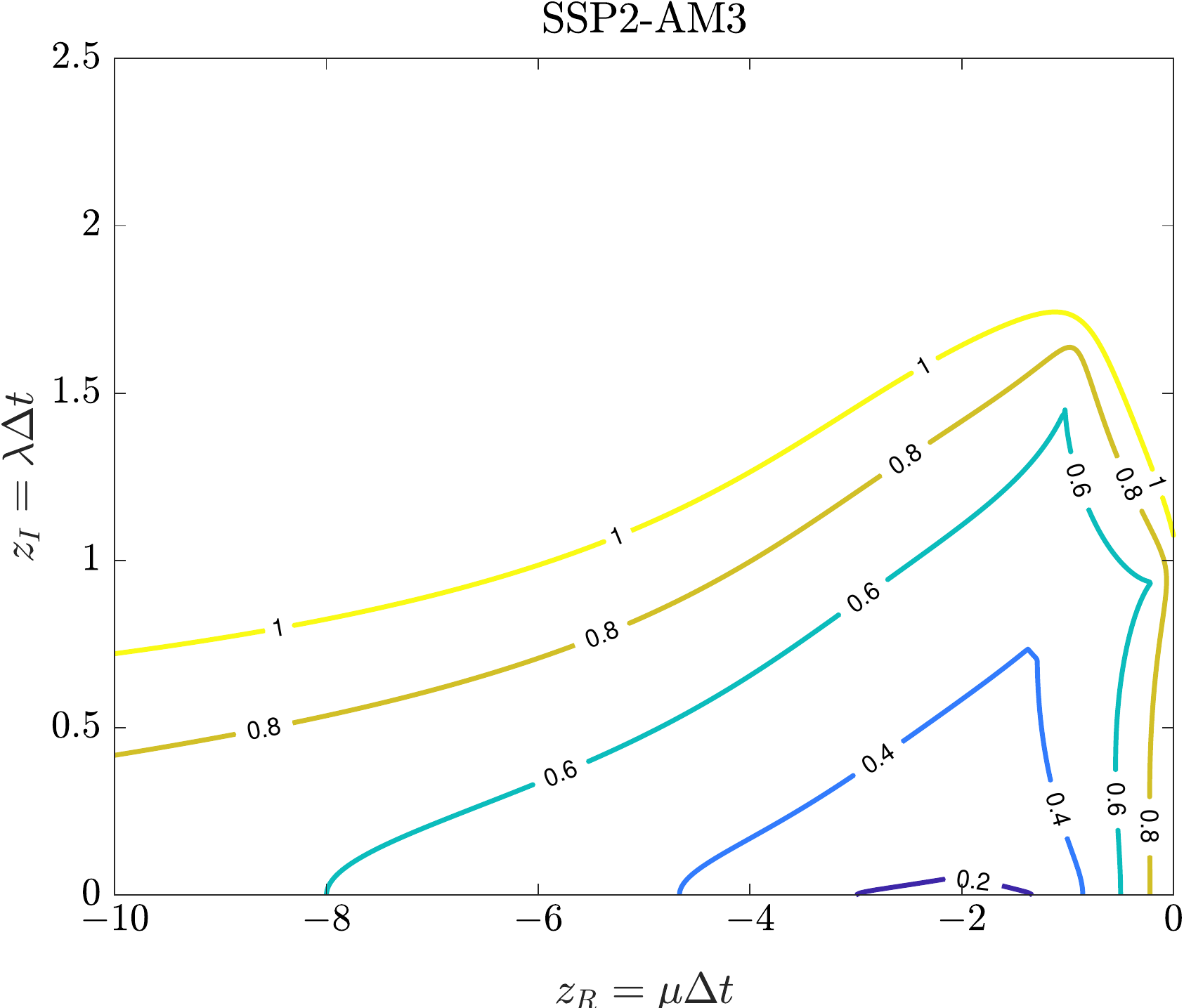}\\[+.1cm]
\includegraphics[scale=0.33]{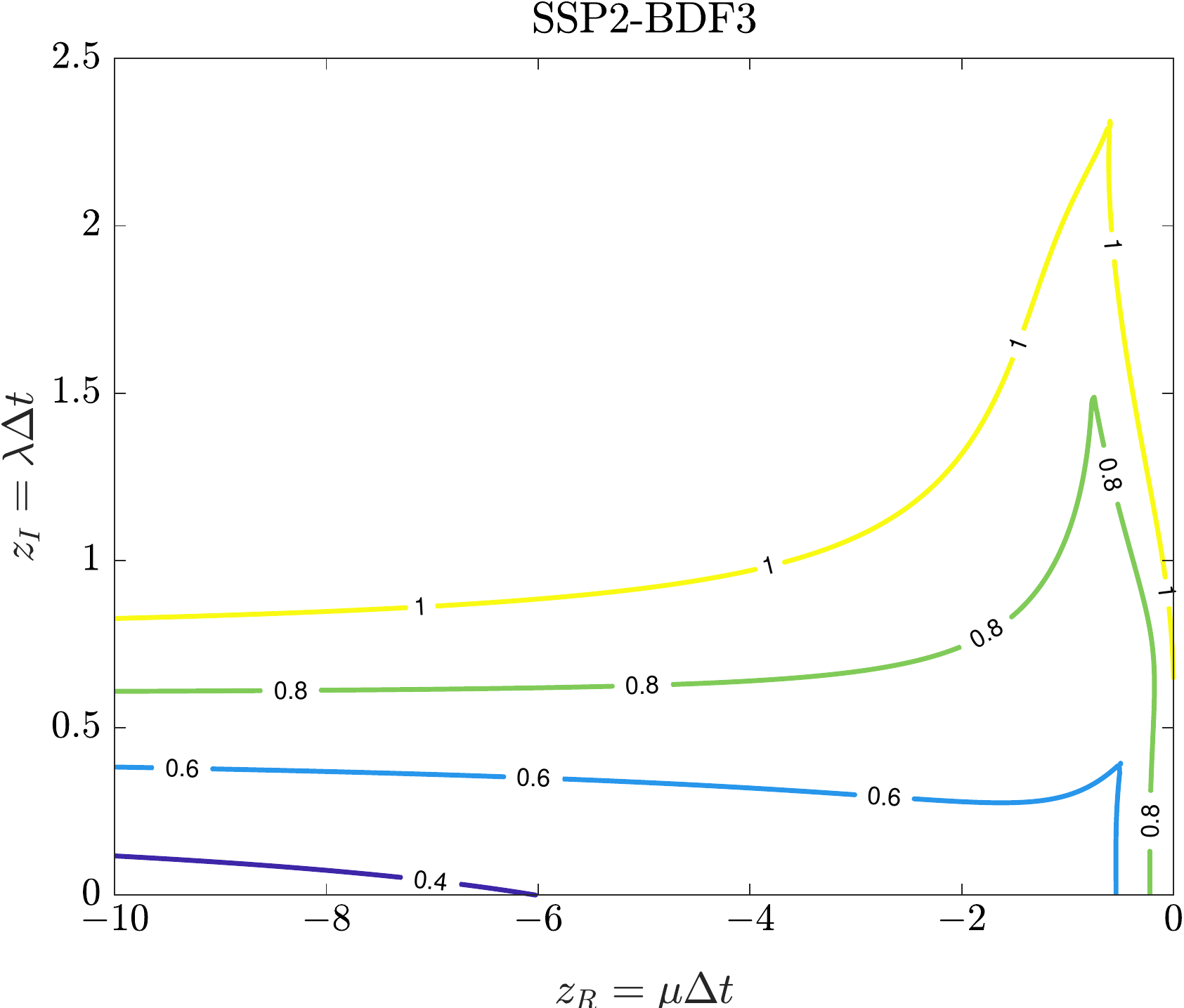}\quad
\includegraphics[scale=0.33]{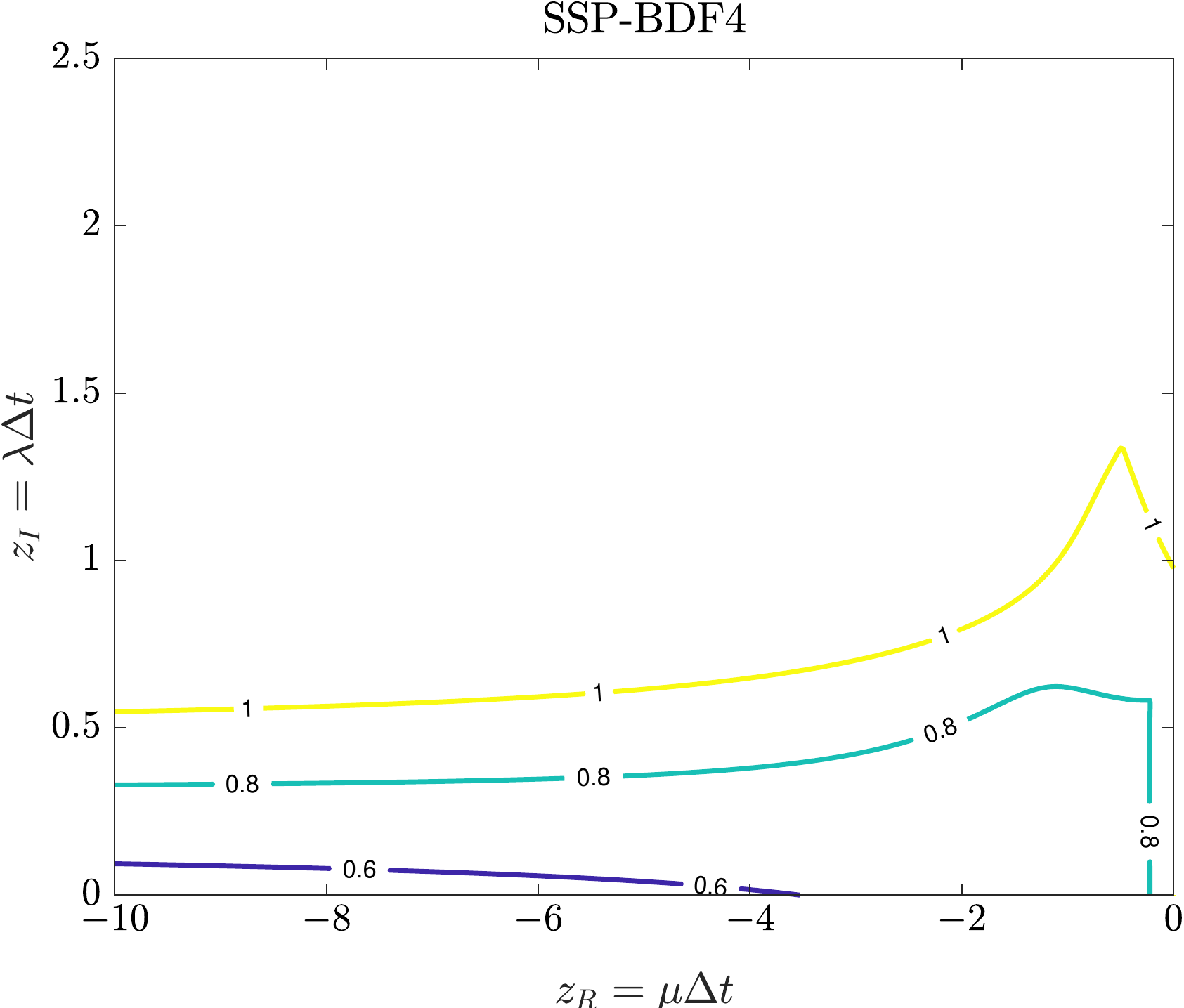}
\end{center}
\caption{Contours of the stability region in the case of problem (\ref{eq:add}) for third order and fourth order SSP methods. The use of the SSP property permits to improve the stability behavior for $\mu \ll 0$.}
\label{fig:SSP1}
\end{figure}

\paragraph{Third order methods and higher.}
The most natural way to obtain third order methods is to combine as explicit solver a two-step Adams-Bashforth method with a two-step Adams-Moulton method. This same strategy actually can be used to obtain methods of higher order. We will refer to this general class of schemes as AB-AM$p$, where $p$ is the order of the resulting scheme. Except for AB-AM2, which is the same as FE-CN2, these schemes in general suffer of poor stability properties when $\mu \ll 0$ (see Figure \ref{fig:SSP1}). Replacing the Adams-Moulton methods with BDF schemes with the same order yields a class of schemes with better stability properties referred to as AB-BDF$p$, where $p$ is the order of the resulting scheme. In this way, AB-BDF2 is the same as FE-BDF2.

\begin{figure}[htb]
\begin{center}
\includegraphics[scale=0.33]{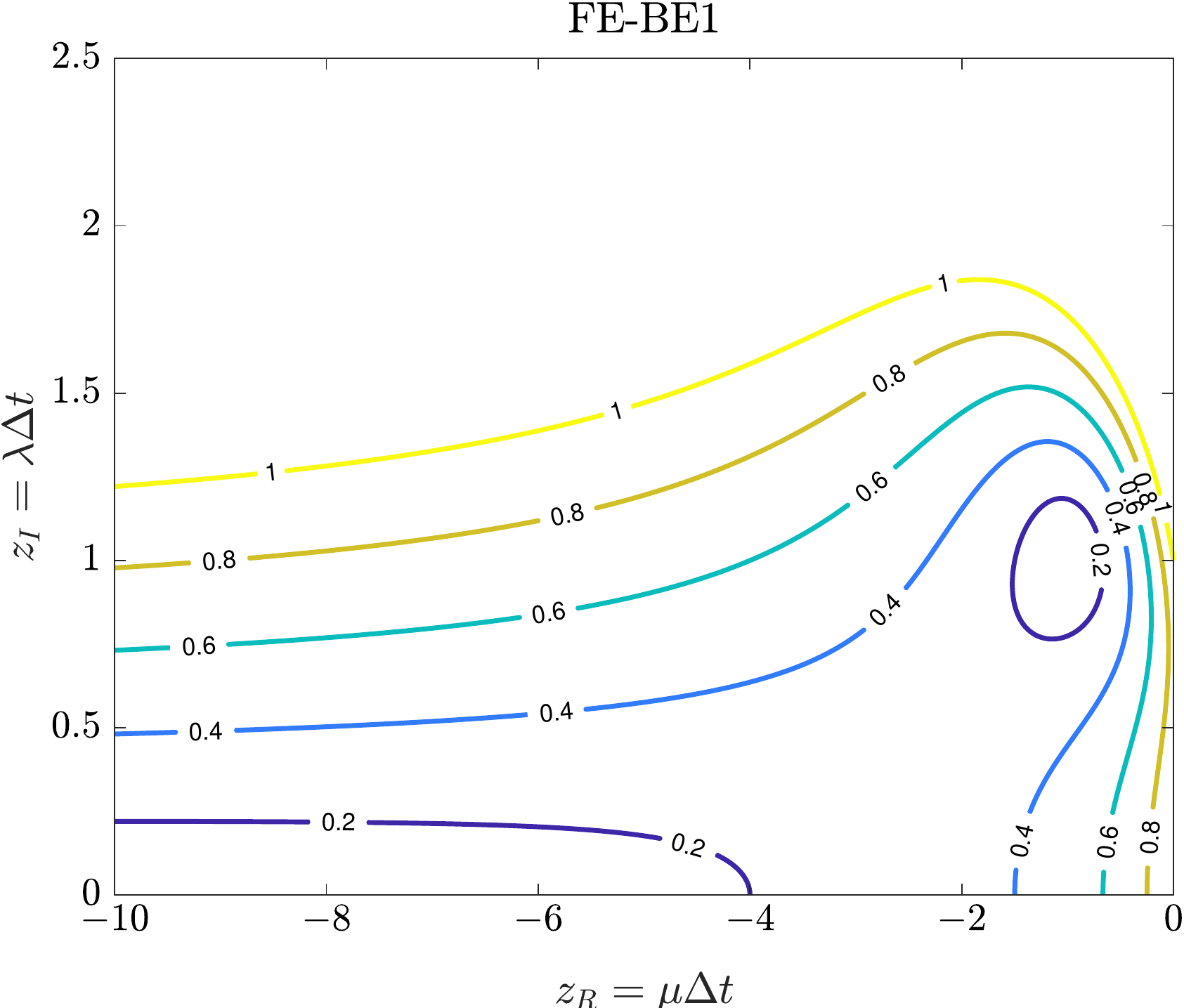}\quad
\includegraphics[scale=0.33]{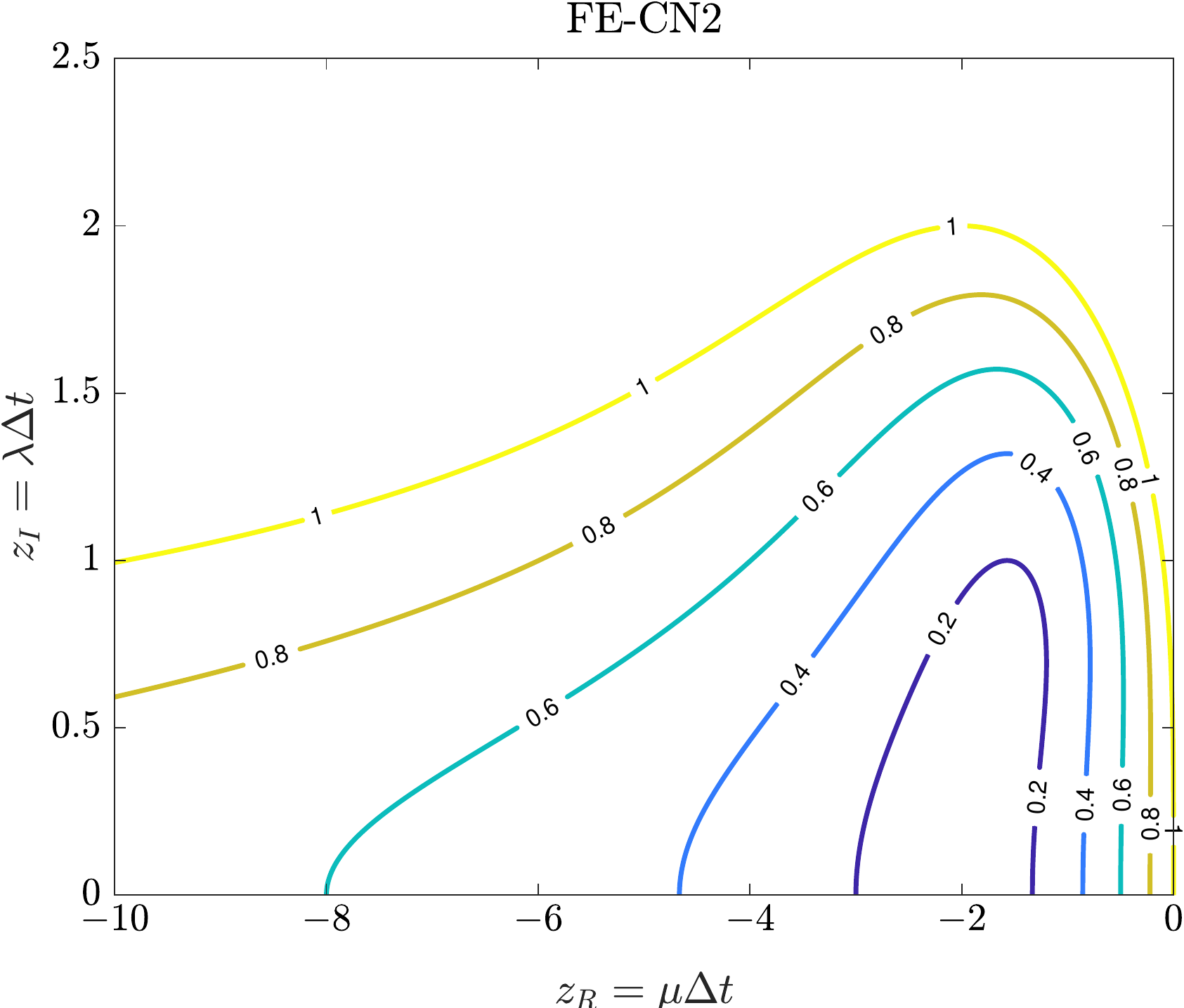}\\[+.1cm]
\includegraphics[scale=0.33]{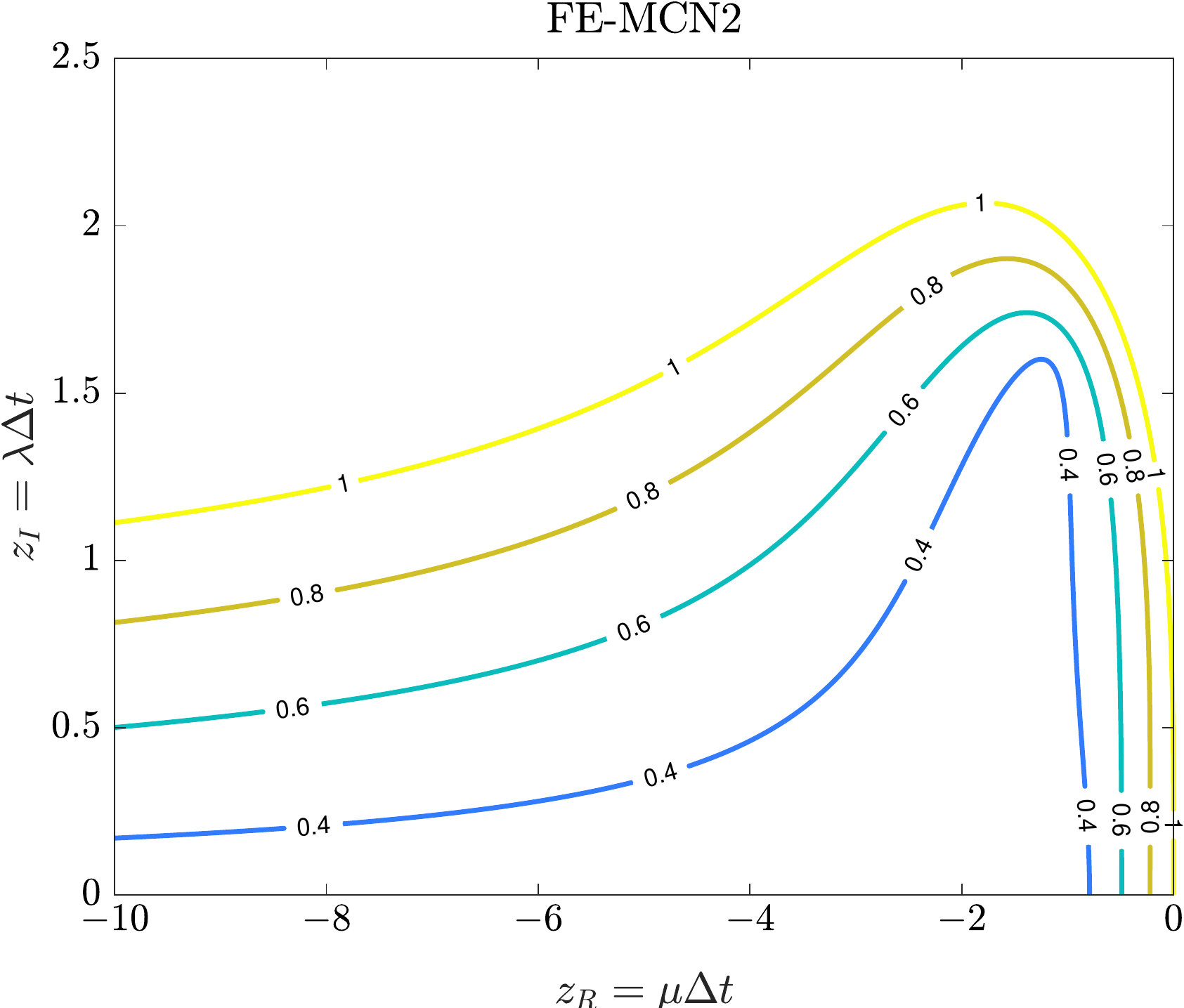}\quad
\includegraphics[scale=0.33]{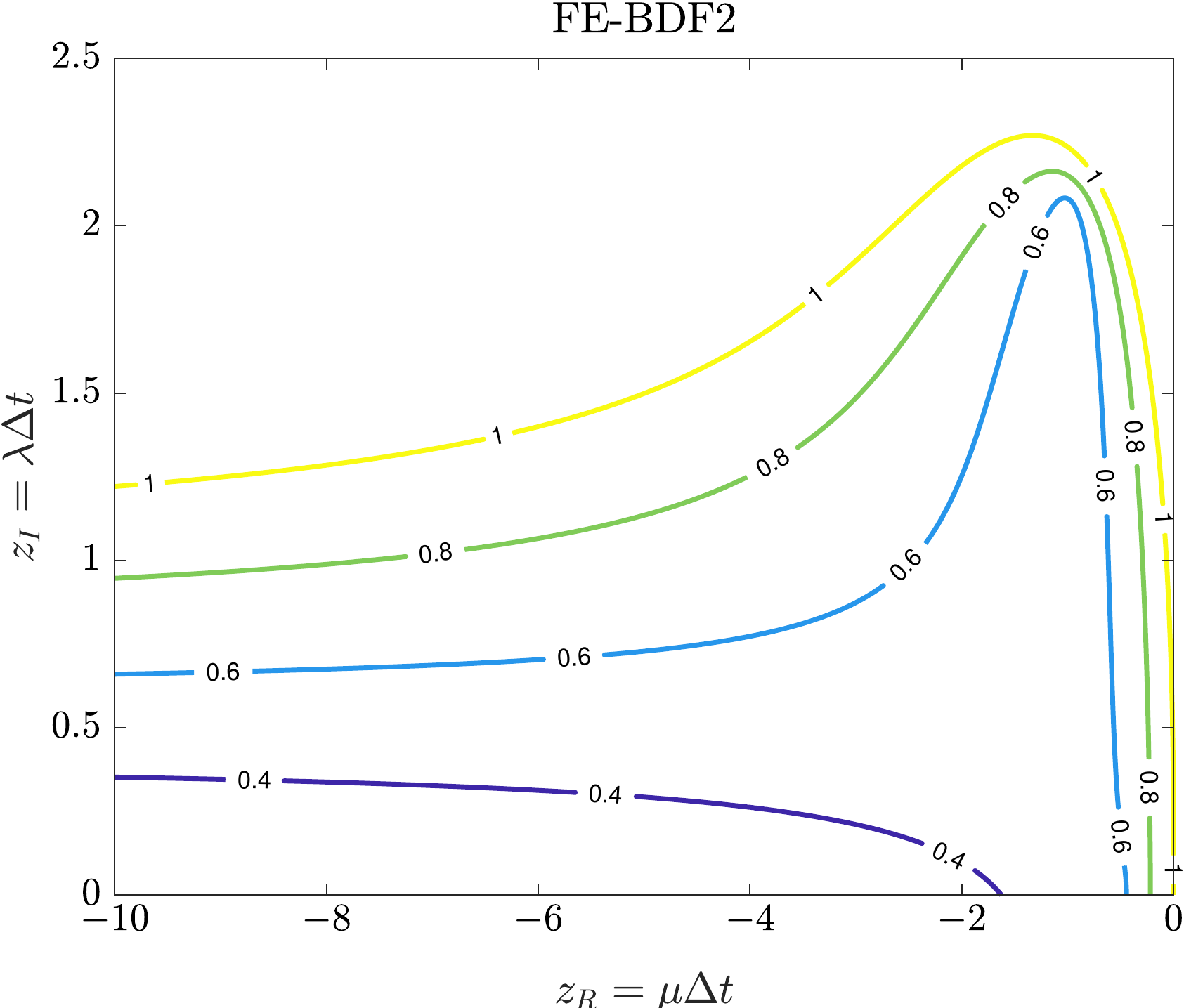}\\[+.1cm]
\includegraphics[scale=0.33]{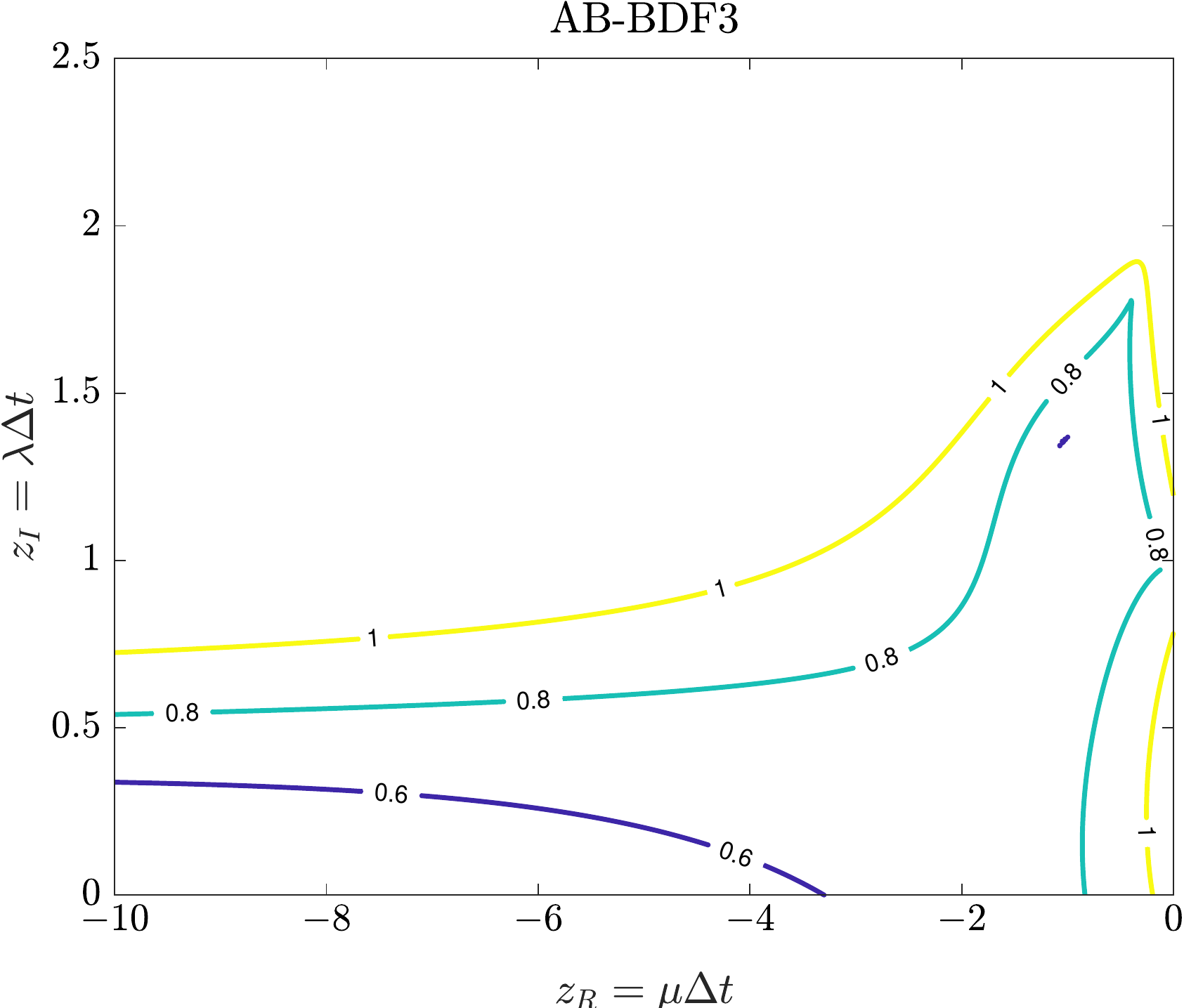}\quad
\includegraphics[scale=0.33]{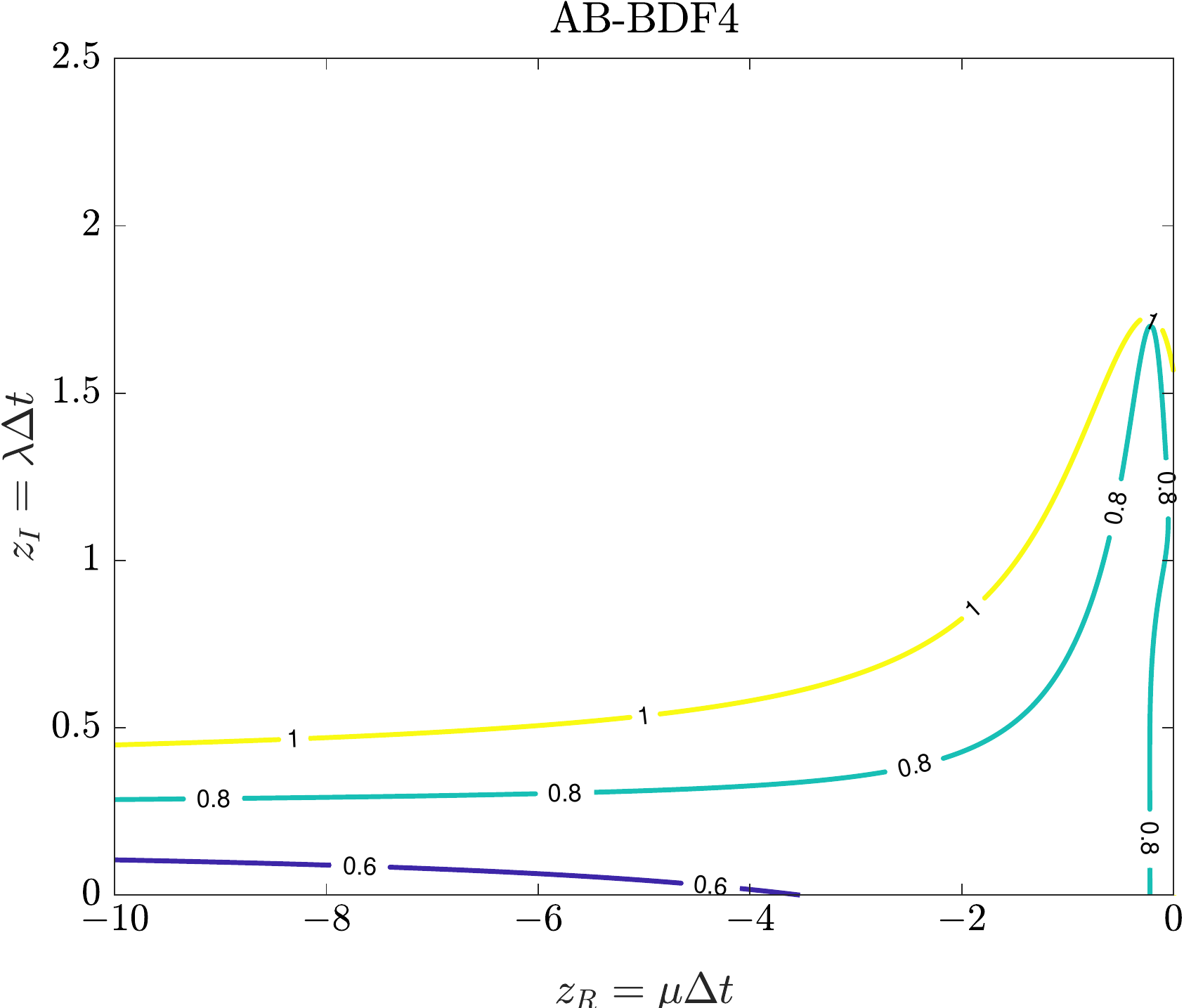}
\end{center}
\caption{Contours of the stability region in the case of problem (\ref{eq:add}) for some semi-implicit multistep methods of order one to four. Increasing the order of accuracy the stability constraints on $\Delta t$ become more severe for $\mu \ll 0$.}
\label{fig:stab1}
\end{figure}

We report in Figure \ref{fig:stab1} the stability contours for various semi-implicit multistep methods up to order four. 

\paragraph{Strong stability preserving methods.} Often the time integration of PDEs requires some monotonicity properties to be satisfied. An important class of methods in this direction is represented by the so-called strong stability preserving (SSP) methods. These methods were designed specifically for solving the ODEs coming from a semi-discrete, spatial discretization
of time dependent partial differential equations, especially hyperbolic PDEs and convection-diffusion problems \cite{GST}. In summary, these schemes are stable for a certain (semi) norm
\be
\|u^{n+1}\| \leq \|u^n\|
\label{eq:SSP}
\ee
under a suitable time step restriction $\Delta t \leq \Delta t_0$. Typically this schemes are applied to the convective part and $\Delta t_0$ refers to the stability constraint, usually referred to as CFL condition, which links $\Delta t_0 = C \Delta t_{FE}$, where $C$ is the CFL coefficient and $\Delta t_{FE}$ the stability constraint of the SSP property in the Forward Euler scheme. 


In \cite{GST} it is shown that there are no implicit multistep SSP schemes of order higher than 1. Therefore, we can combine optimal explicit multistep methods satisfying the SSP property \eqref{eq:SSP} as predictor for the non stiff component with implicit methods to improve the overall stability region in our semi-implicit schemes. 

The explicit multistep SSP schemes are of the general form
\be
u^{n+1} = \sum_{j=0}^{s-1} \left(\alpha_j u^{n-j} + \Delta t \beta_j {\mathcal H}\left({u}^{n-j},v^{n-j}\right)
\right),\qquad \alpha_j \geq 0.
\ee
The optimal second order two steps explicit SSP method ($C=1/2$) corresponds to the choices ${\boldsymbol \alpha}=(4/5,1/5)^T$ and ${\boldsymbol \beta}=(8/5,-2/5)^T$, whereas the optimal third order explicit four steps scheme ($C=1/3$) is obtained for ${\boldsymbol \alpha}=(16/27,0,0,11/27)^T$ and ${\boldsymbol \beta}=(16/9,0,0,4/9)^T$ (see \cite{GST}). The corresponding schemes are denoted as SSP-AM3, SSP-BDF3 and SSP-BDF4.
If one increases the number of steps, then SSP methods can be found to have larger SSP regions. Because there is no significant increase in the computational cost when the number of steps is increased, if storage is not a consideration, it may be advantageous to use a SSP multi-step methods with more steps and larger stability domain. For example, the optimal second order four steps SSP method ($C=2/3$) is obtained with ${\boldsymbol \alpha}=(8/9,0,0,1/9)^T$ and ${\boldsymbol \beta}=(4/3,0,0,0)^T$. The corresponding third order semi-implicit schemes are denoted as SSP2-AM3 and SSP2-BDF3.
 
We report in Figure \ref{fig:SSP1}, the stability contours of AB-AM$3$ and SSP-AM$3$. The importance of the optimal SSP property in the explicit scheme is evident and improves dramatically the stability properties of the resulting method. For large values of $|\mu|$ the semi-implicit method based on AM3 becomes stable under a reasonable CFL condition of $\Delta t z_I \leq 1/2$. Similarly the use of SSP2 predictor increase this stability constraint to approximatively $\Delta t z_I \leq 3/5$. We also gave in the same Figure the optimal third and fourth order methods constructed using the BDF formulae. Again the use of SSP2 predictor slightly improves the stability properties for large $|\mu|$.






\section{Numerical results}
In what follows we investigate semi-implicit schemes for non-linear problems with stiff terms, where the stiffness can be solve efficiently by linear solvers. 

We will discuss non-linear diffusion-reaction and convective-diffusion problems, following different examples extracted from \cite{BFR} and \cite{HV}.

\subsection{Test 1: Order of convergence for reaction-diffusion system}
Following the validation example of \cite{BFR} we consider the non-autonomous diffusion-reaction system, where $\omega= (\omega_1,\omega_2):\mathbb{R}_+\times[0,2\pi)^2\to\mathbb{R}^2$ is solution of the following system
\begin{equation}\begin{aligned}\label{eq:GrayScott_mod}
\partial_t \omega_1 &= \Delta \omega_1 -\alpha(t)\omega_1^2 +\frac{9}{2}\omega_1+\omega_2+f(t)\\
\partial_t \omega_2 &= \Delta \omega_2 +\frac{7}{2}\omega_2,
\end{aligned}\end{equation}
the time dependent factors are $\alpha(t)=2 e^{t/2}$ and $f(t) = -2 e^{-t/2}$. 
Accounting periodic boundary conditions, the initial data is extracted from the exact solution
\begin{equation}\begin{aligned}\label{eq:GrayScott_mod0}
\omega_1(t,x,y)&=  e^{-t/2}\left(1+\cos(x)\right)\\
\omega_2(t,x,y)&= e^{-t/2}\cos(2x)
\end{aligned}
\qquad (x,y)\in[0,2\pi)^2
\end{equation}

In order to apply the semi-implicit multi-step scheme \eqref{eq:semiGPS} we reformulate system \eqref{eq:GrayScott_mod} introducing $u=(u_1,u_2)$ and $v=(v_1,v_2)$, and the operator
\[
\mathcal{H}(t,u,v) =
\begin{pmatrix}
\Delta v_1 -\alpha(t)u_1v_1 +\frac{9}{2}v_1+v_2+f(t)\cr~\cr
\Delta v_2 +\frac{5}{2}v_2
\end{pmatrix}.
\]


For the spatial discretization of the diffusion operator we apply a 6th order central finite difference, on an uniform grid for the periodic domain $[0,2\pi)^2$ with $\Delta x = \Delta y$. For $u(t,x,y)$ evaluated on a point $(t_n,x_i,y_j)$ we consider the operators
\begin{equation*}\begin{aligned}
(D_{xx} u^n)_{ij}= \frac{2u^n_{i+3j}-27u^n_{i+2j}+270u^n_{i+1j}-240u^n_{ij}+270u^n_{i-1j}-27u^n_{i-2j}+2u^n_{i-3j}}{180\Delta x^2}.
\end{aligned}\end{equation*}
 and similarly $D_{yy}$ for the $y-direction$, taking into account periodic boundary conditions.

In order to estimate the order of accuracy in time we proceed refining the space step and time step with stability conditions
\[
\Delta t = \lambda{\Delta x}, 
\]
where we choose $\lambda = 0.5$. Thus we monitor the $\ell_1$ error decay of the numerical solution $\omega_{ij}(t)$ at final time $T = 2$
\[
\ell_1(\omega^{(k)}) = \Delta x\Delta y\sum_{i,j}\|\omega^{(k)}_{ij}(T)-\omega(x_i,y_j,T)\|,
\]
where $\omega^{(k)}_{ij}(T)$ is the numerical solution evaluated on $N_x=2^k$ points, with $k = 6,7,8$ and $\omega(x_i,y_j,t_n)$ is the exact solution evaluated at $(x_i,y_j,T)$.

%

\begin{figure}
	\centering
	\subfigure[$\omega=(\omega_1,\omega_2)$]{\includegraphics[width=5.5cm]{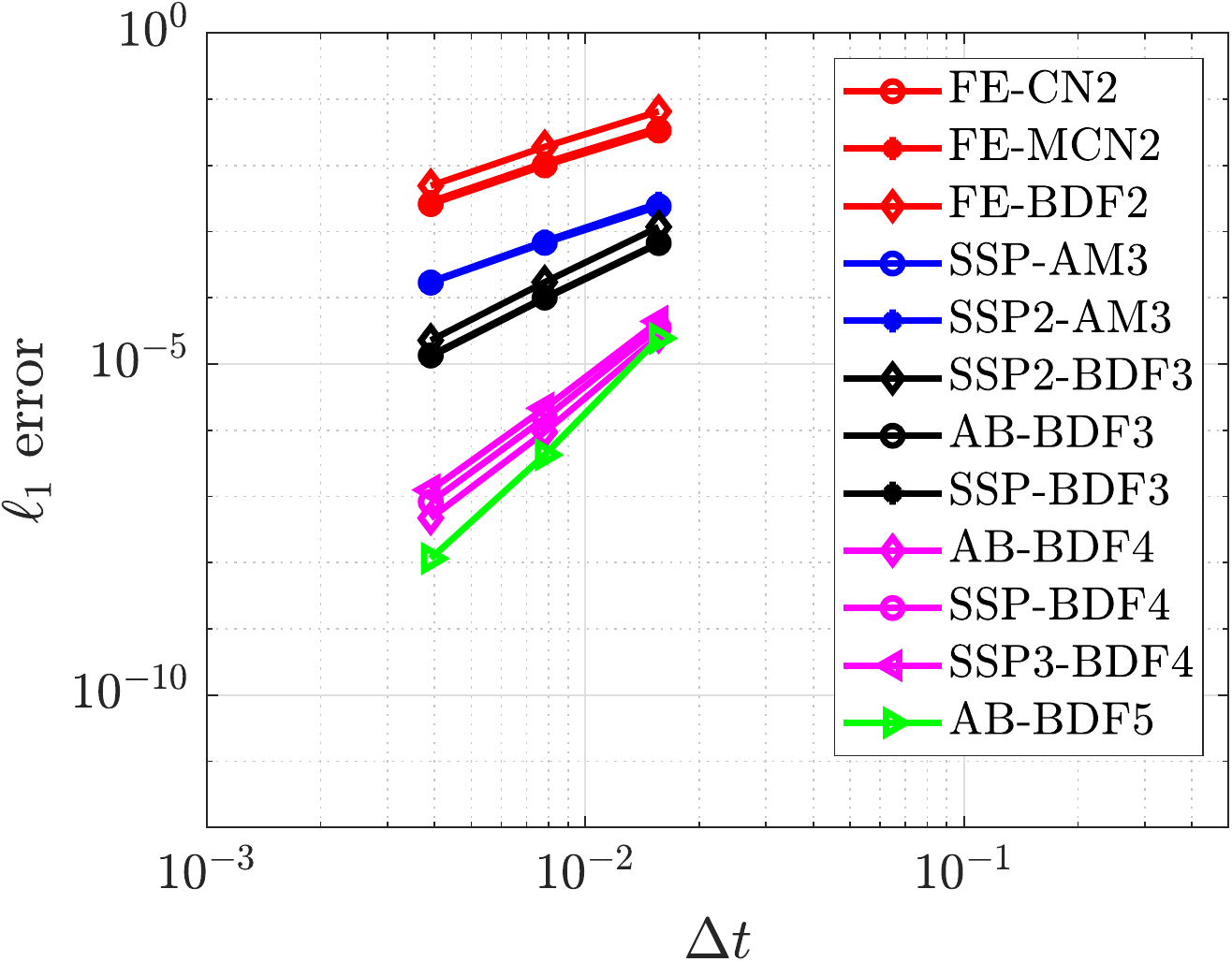}}\hspace{0.5cm}
	\subfigure[$\omega_2$]{\includegraphics[width=5.5cm]{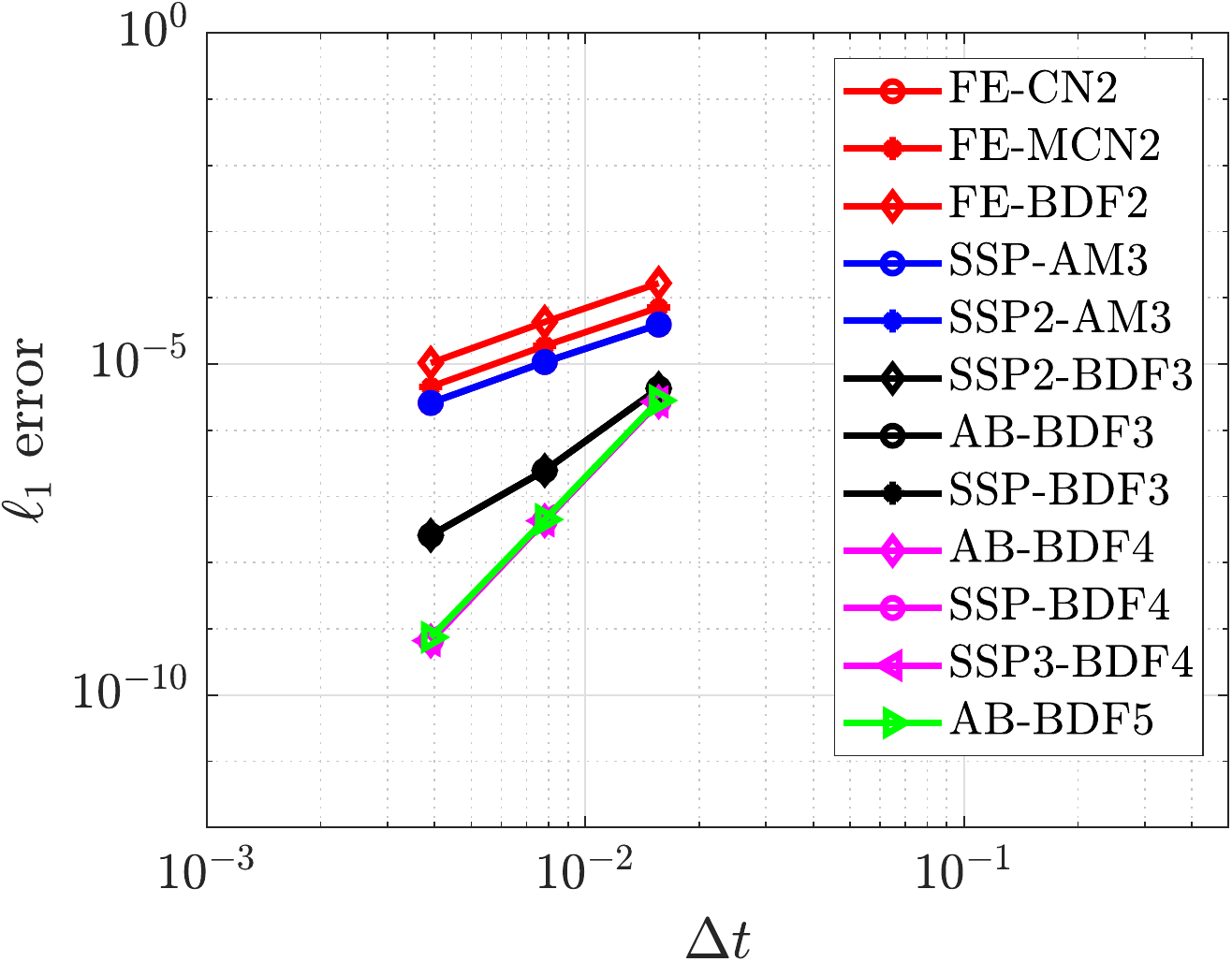}}
	\caption{Test 1: We report $\ell_1$ error norm of the solution of \eqref{eq:GrayScott_mod} obtained for different schemes from 2nd to 5th order. The right plot (a) depicts the order of convergence for the full solution $\omega$. Plot (b) shows the error relative to the autonomous component $\omega_2$.}\label{fig:errorl1}
\end{figure}

In Figure \ref{fig:errorl1} we report the $\ell_1$ error decay for several schemes. We compare the accuracy of the solution (a) $\omega=(\omega_1,\omega_2)$, on the left side, with respect to (b) the accuracy of $\omega_2$, observing better convergence properties for the second component which non-autonomous and independent of $\omega_1$. 
Nonetheless in both cases 5th order of convergence is observed.

\subsection{Test 2: Non-linear reaction diffusion system}
We consider the Gray-Scott model, as studied in \cite{Pe,ZWZ,HV}, 
\begin{equation}\begin{aligned}\label{eq:GrayScott}
\partial_t \omega_1 = \sigma_1\Delta \omega_1 -\omega_1\omega_2^2 +\gamma(1-\omega_1)\\
\partial_t \omega_2 = \sigma_2\Delta \omega_2 +\omega_1\omega_2^2 -(\gamma+\kappa)\omega_2
\end{aligned}\end{equation}
where $\omega= (\omega_1,\omega_2):\mathbb{R}_+\times[-1,1)^2\to\mathbb{R}^2$ with periodic boundary conditions, and initial data
\begin{equation}\begin{aligned}\label{eq:GrayScott_v0}
\omega_1(0,x,y)&= 1-2\omega_2(0,x,y)\\
\omega_2(0,x,y)&= \frac{1}{4}\sin^2(4\pi x)\sin^2(4\pi y)
\end{aligned}
\qquad (x,y)\in[-1/4,1/4]^2
\end{equation}
and zero otherwise on the squared domain $[-1,1]^2$. The diffusion and reaction parameters are 
\[\sigma_1 = 8\times10^{-5},\quad \sigma_2 = 4\times10^{-5},\quad \gamma = 0.024,\quad \kappa = 0.06.\]
Initial data and parameters are selected in order to match the test proposed in \cite{HV} and inspired from \cite{Pe}.
In order to apply the semi-implicit multi-step scheme \eqref{eq:semiGPS} we reformulate system \eqref{eq:GrayScott} introducing $u=(u_1,u_2)$ and $v=(v_1,v_2)$, and the operator
\[
\mathcal{H}(t,u,v) =
\begin{pmatrix}
D_1\Delta u_1 -v_1u_2^2 +\gamma(1-v_1)\cr
D_2\Delta u_2 +v_1u_2^2 -(\gamma+\kappa)v_2
\end{pmatrix},
\]
where now we treat explicitly the diffusion terms, since $D_1,D_2$ are negligible with respect to the reaction coefficients.

Thus the linear reaction terms are taken implicitly, whereas the non-linear term is taken implicit only in the $\omega_1$ component. We employ LM scheme SSP3-BDF4 for time integration and forth order central difference for the space discretization, introducing the operators
\begin{align}\label{eq:diff4}
(D_{xx} u^n)_{ij} &= \frac{-u^n_{i+2j}+16u^n_{i+1j}-30u^n_{ij}+16u^n_{i-1j}-u^n_{i-2j}}{12\Delta x^2}.
\end{align}
 and similarly $D_{yy}$ for the $y-direction$, taking into account periodic boundary conditions.

\begin{figure}
	\centering
	\includegraphics[width=4.95cm]{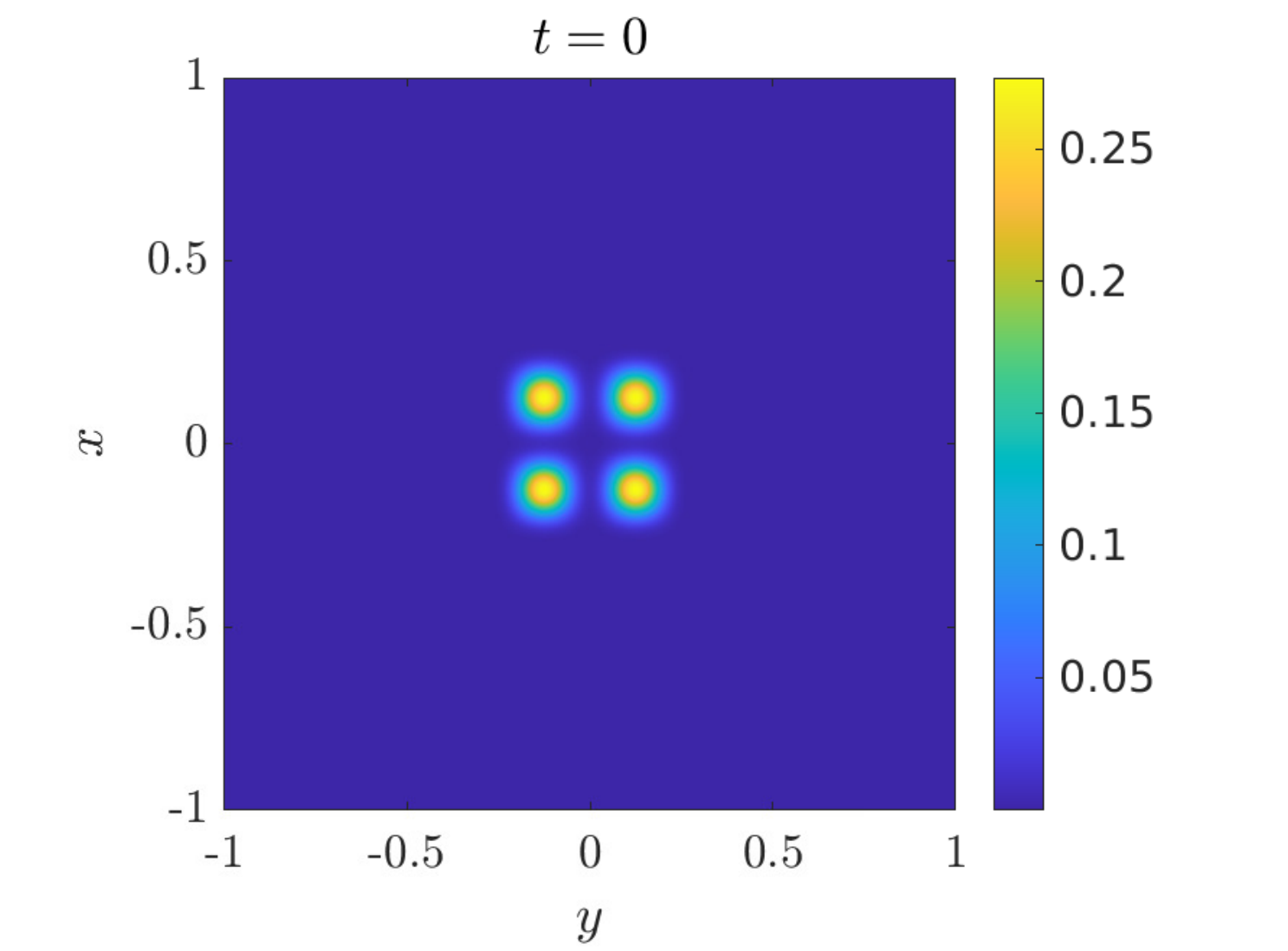}	\includegraphics[width=4.95cm]{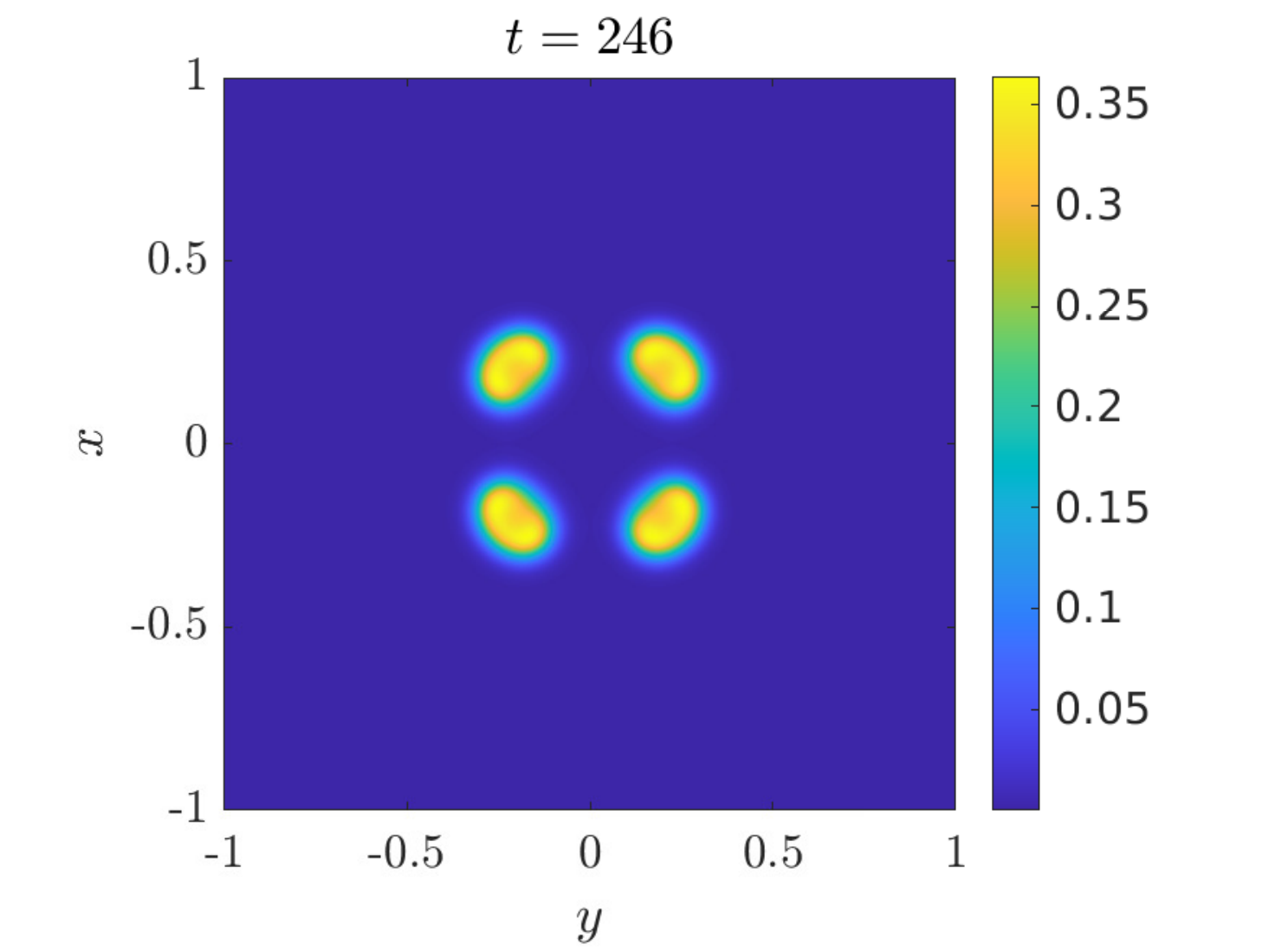} \includegraphics[width=4.95cm]{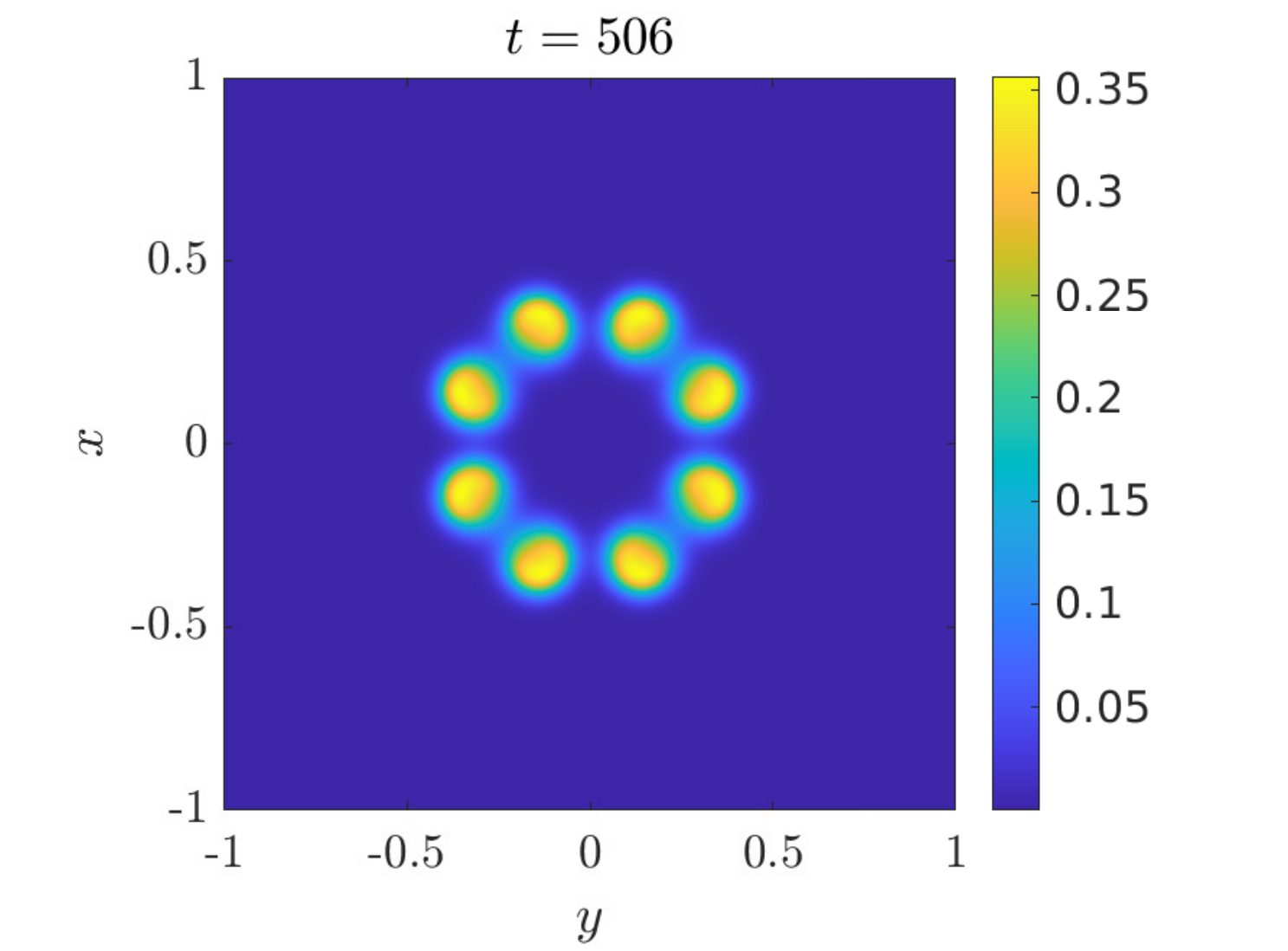}
	\\
	\includegraphics[width=4.95cm]{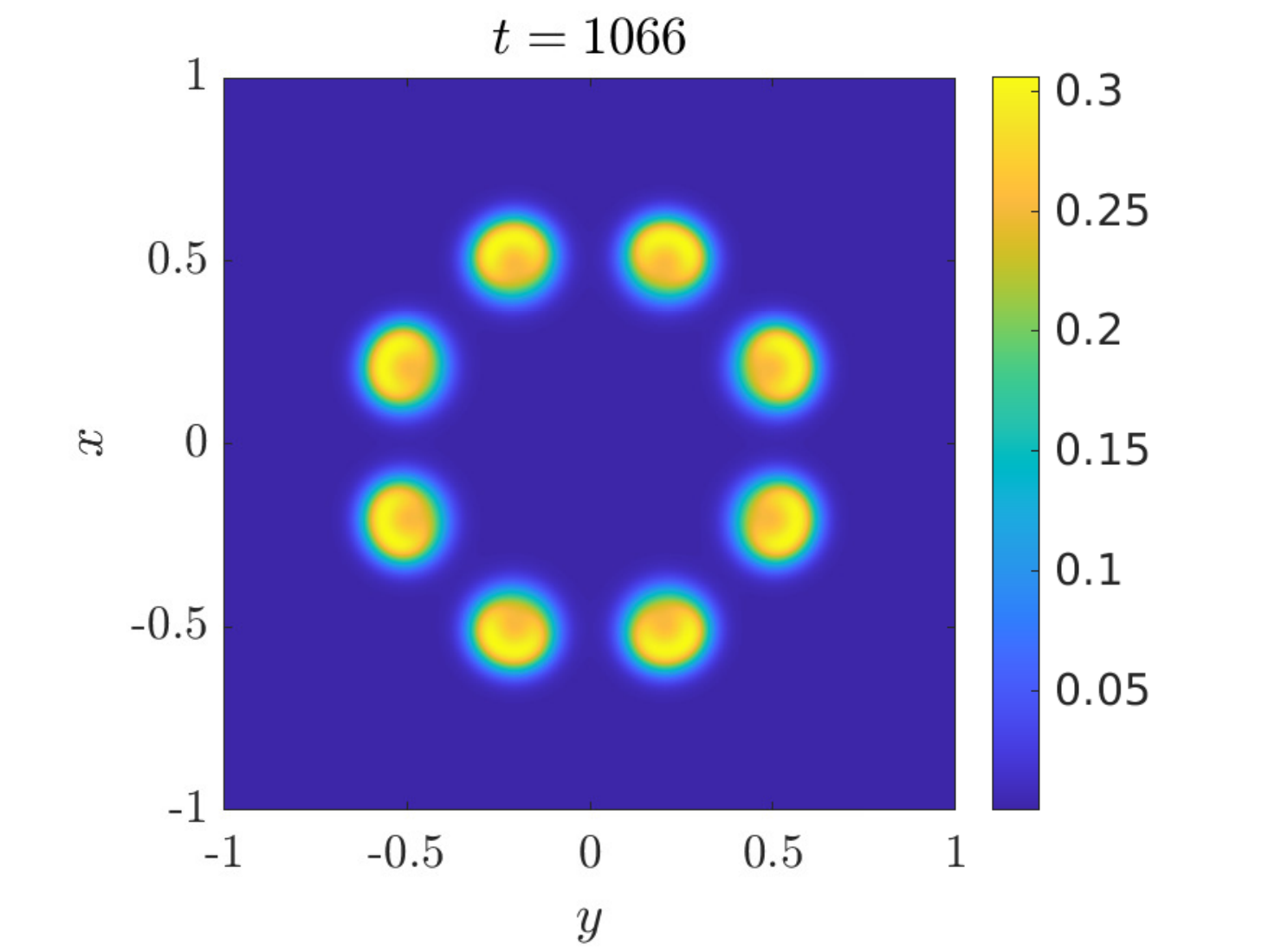}	\includegraphics[width=4.95cm]{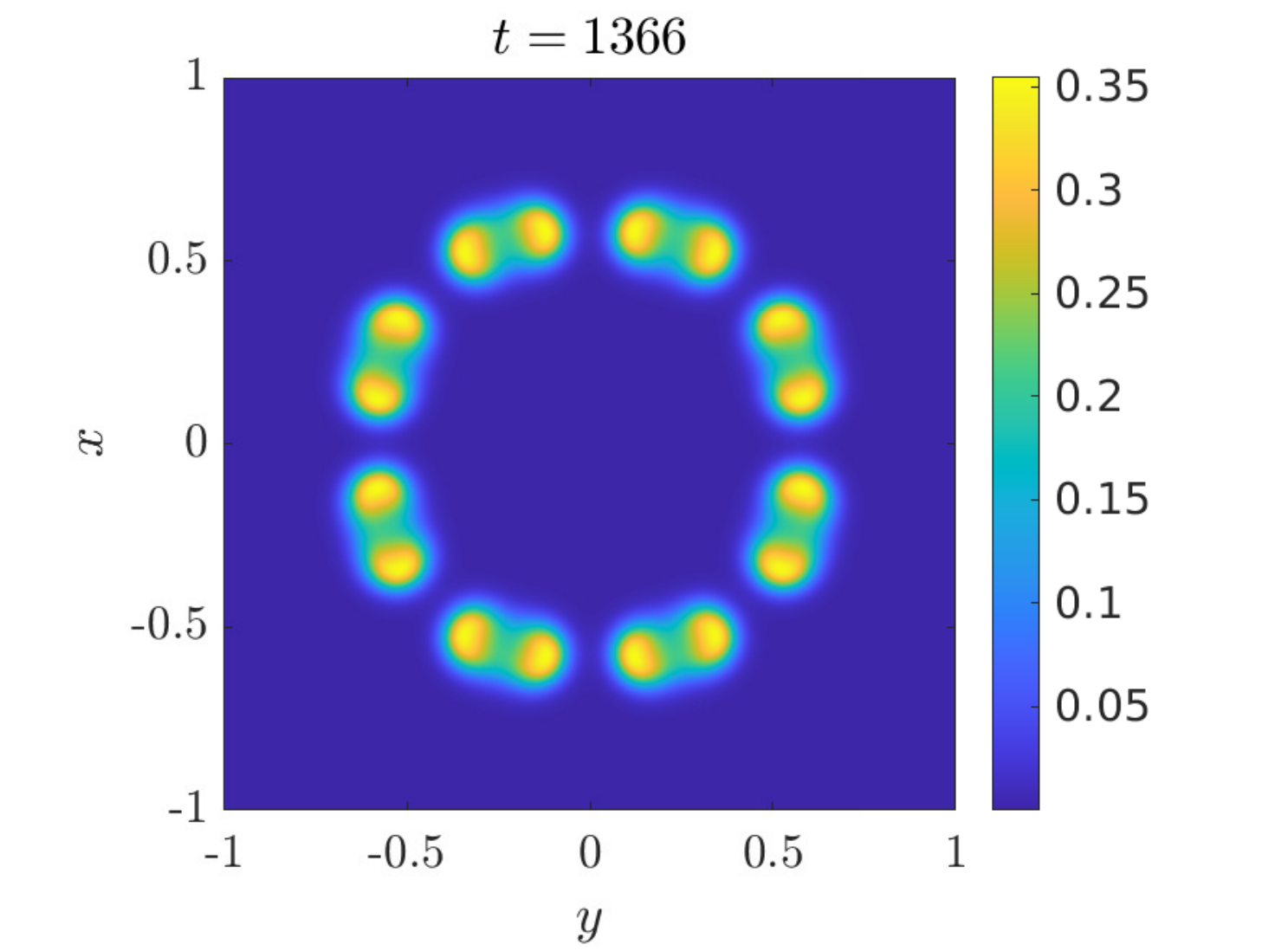} \includegraphics[width=4.95cm]{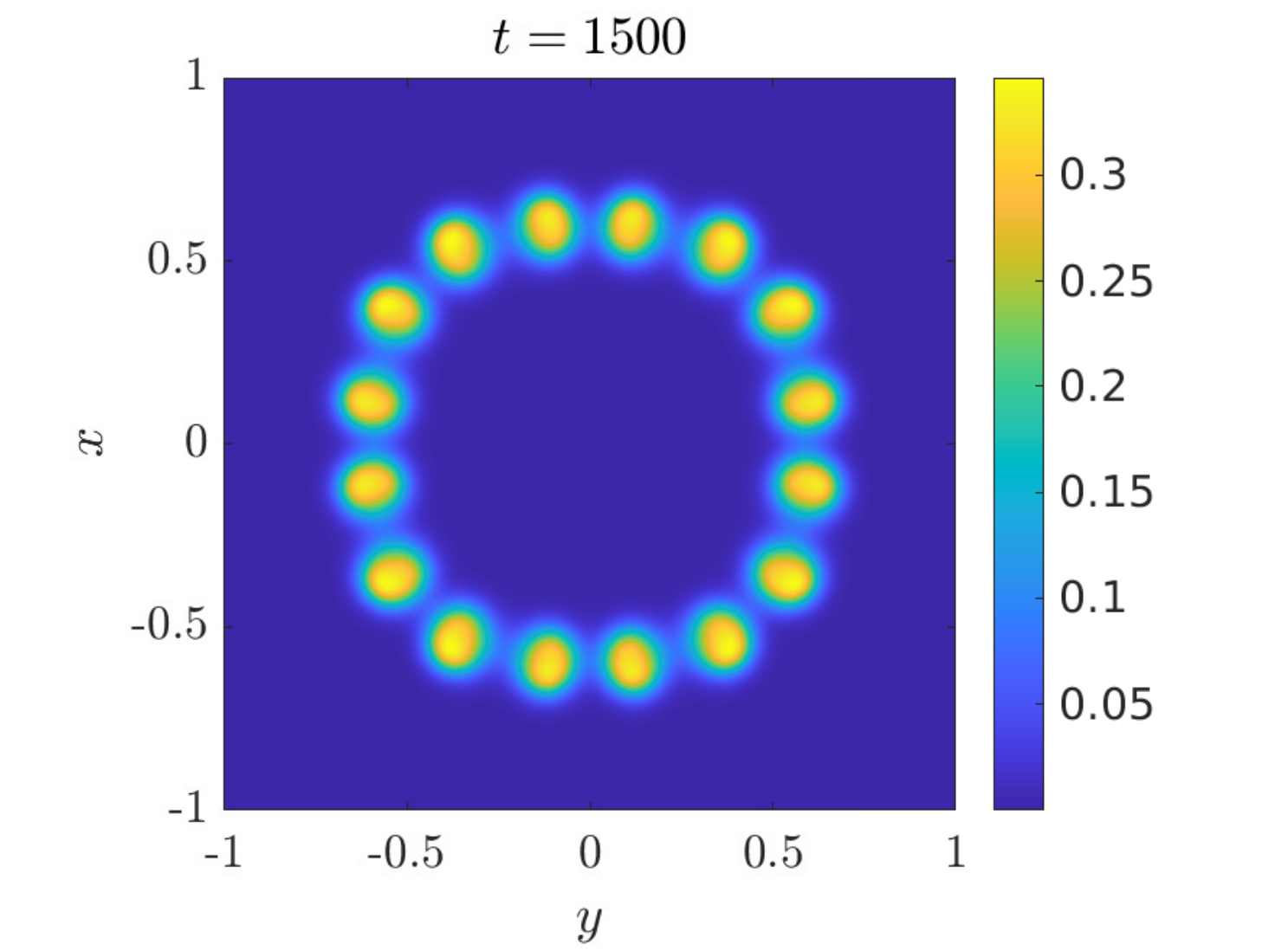}
	\caption{Test 2: {(\em Gray-Scott model).} Solution of the reaction-diffusion system \eqref{eq:GrayScott} in the component $\omega_2$ at different time frames obtained with scheme \textit{SSP-BDF4} on the square $\Omega=[-1, 1]^2$ with uniform space-time discretization with $\Delta x=0.01$ and $\Delta t=\Delta x$, and periodic boundary conditions. Parameters of the model are  $\sigma_1=8\times10^{-5}$, $\sigma_2=4\times10^{-5}$, $\gamma=0.024$ and $\kappa=0.06$.} \label{fig:GS_circ}
\end{figure}

Figure \ref{fig:GS_circ} reports the evolution of the  component $\omega_2(t,x,y)$ for different times. Starting from a symmetric concentration Gray-Scott model produces spot multiplication, resembling cell division process. The computational domain is discretized with $N_x=N_y=200$ points in space, final time is $T=1500$ with uniform step $\Delta t = \Delta x/2$.

\newpage

\subsection{Test 3: Nonlinear convection-diffusion}
We finally consider the nonlinear convection diffusion model defined on the full plan $x\in\mathbb{R}^2$ as follows
\begin{equation}\begin{aligned}\label{eq:AD}
&\partial_t \omega +\left(E+\mu\nabla\log(\omega)\right)\cdot\nabla \omega = \mu\Delta \omega\cr
&\omega_0(0,x) = e^{-\|x\|^2/2}
\end{aligned}\end{equation}
where $E=(1,1)^\top$ and $\mu=0.5$. The initial data is extracted from the exact solution given by
\[
\omega(t,x) = \frac{1}{\sqrt{4\mu t+1}}e^{-\frac{\|x- E t\|^2}{8\mu t+2}}.
\]
In order to solve numerically \eqref{eq:AD} we introduce the operator
\[
\mathcal{H}(t,u,v) = -\left(E+\mu\nabla\log(u)\right)\cdot\nabla v + \mu\Delta v
\]
where we treat the convection and diffusion terms implicitly, on the computational domain
$[-10,10]$ up to final time $T=1$. We choose uniform space step $\Delta x=\Delta y= 0.1$ and $\Delta t =\Delta x/2$  For time integration we employ SSP3-BDF4, and 4th order central difference for diffusion operator, as in \eqref{eq:diff4}, and for the convective term we account the operator
\begin{align*}
(D_x u^n)_{ij} &= \frac{-u_{i+2j}+8u_{i+1j}-8u_{i-1j}+u_{i-2j}}{12\Delta x},
\end{align*}
and equivalently $D_y$ for the $y-$direction.

\begin{figure}
	\centering
	\includegraphics[width=4.95cm]{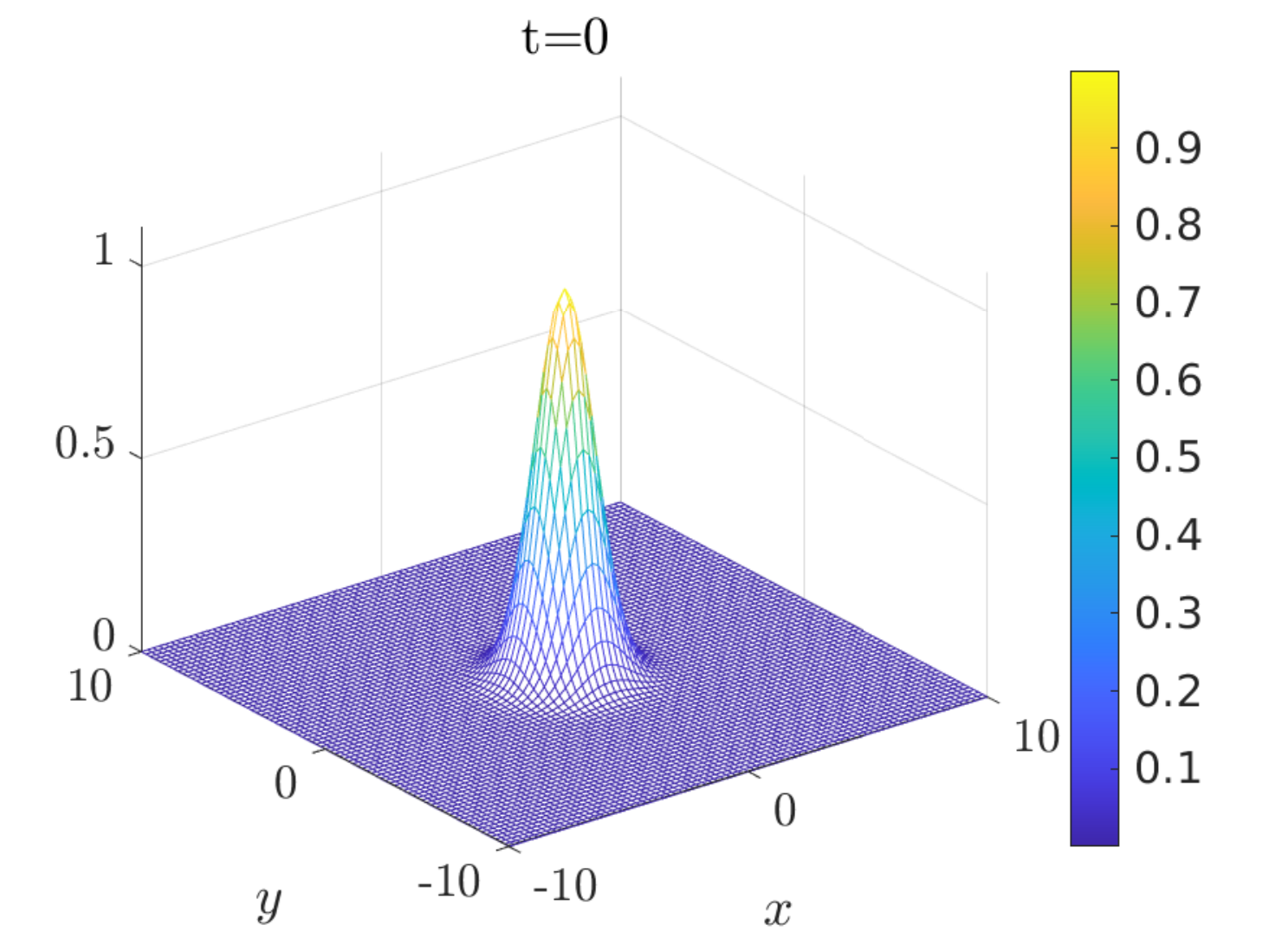}	
	\includegraphics[width=4.95cm]{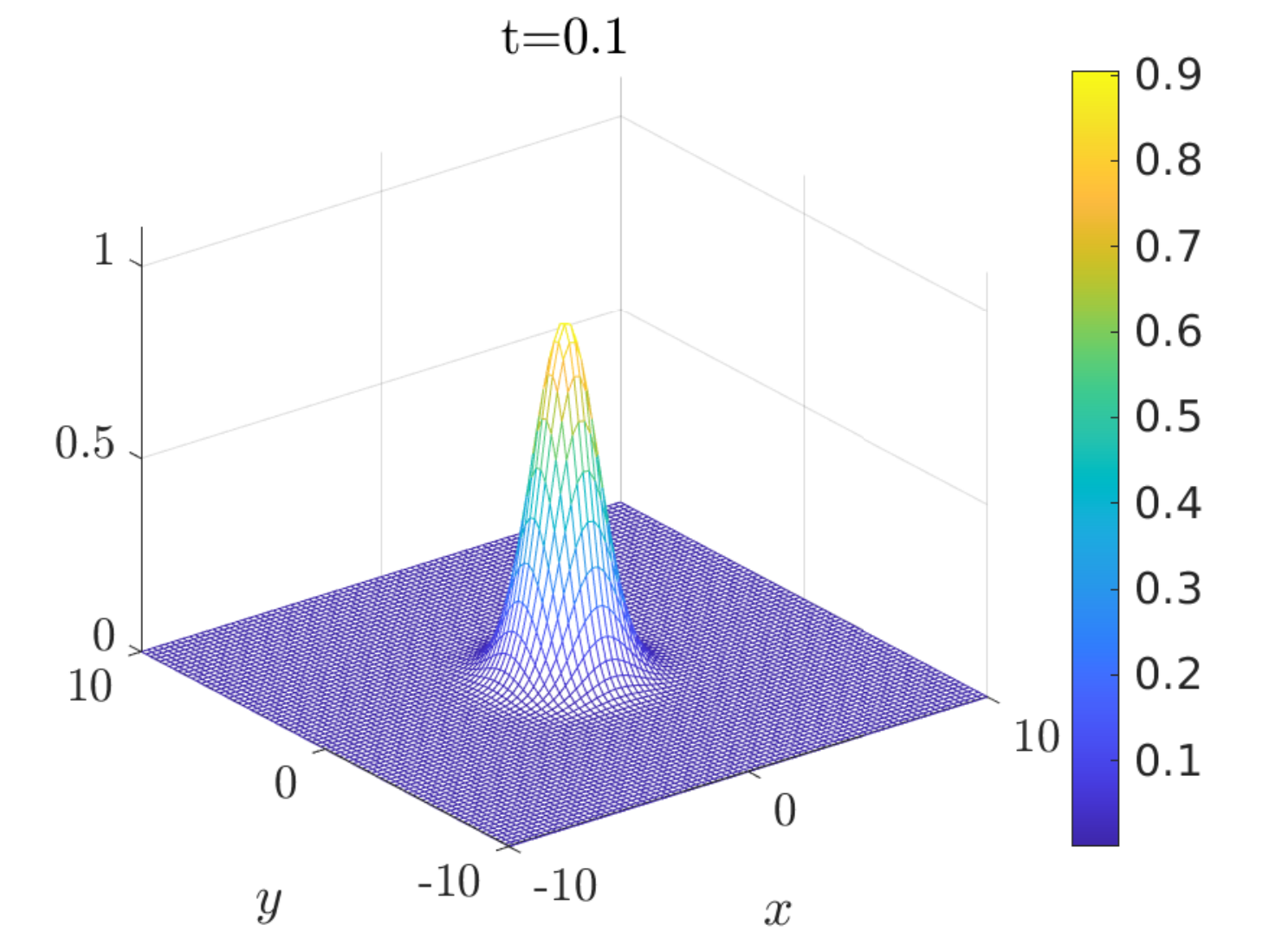}
	\\
	\includegraphics[width=4.95cm]{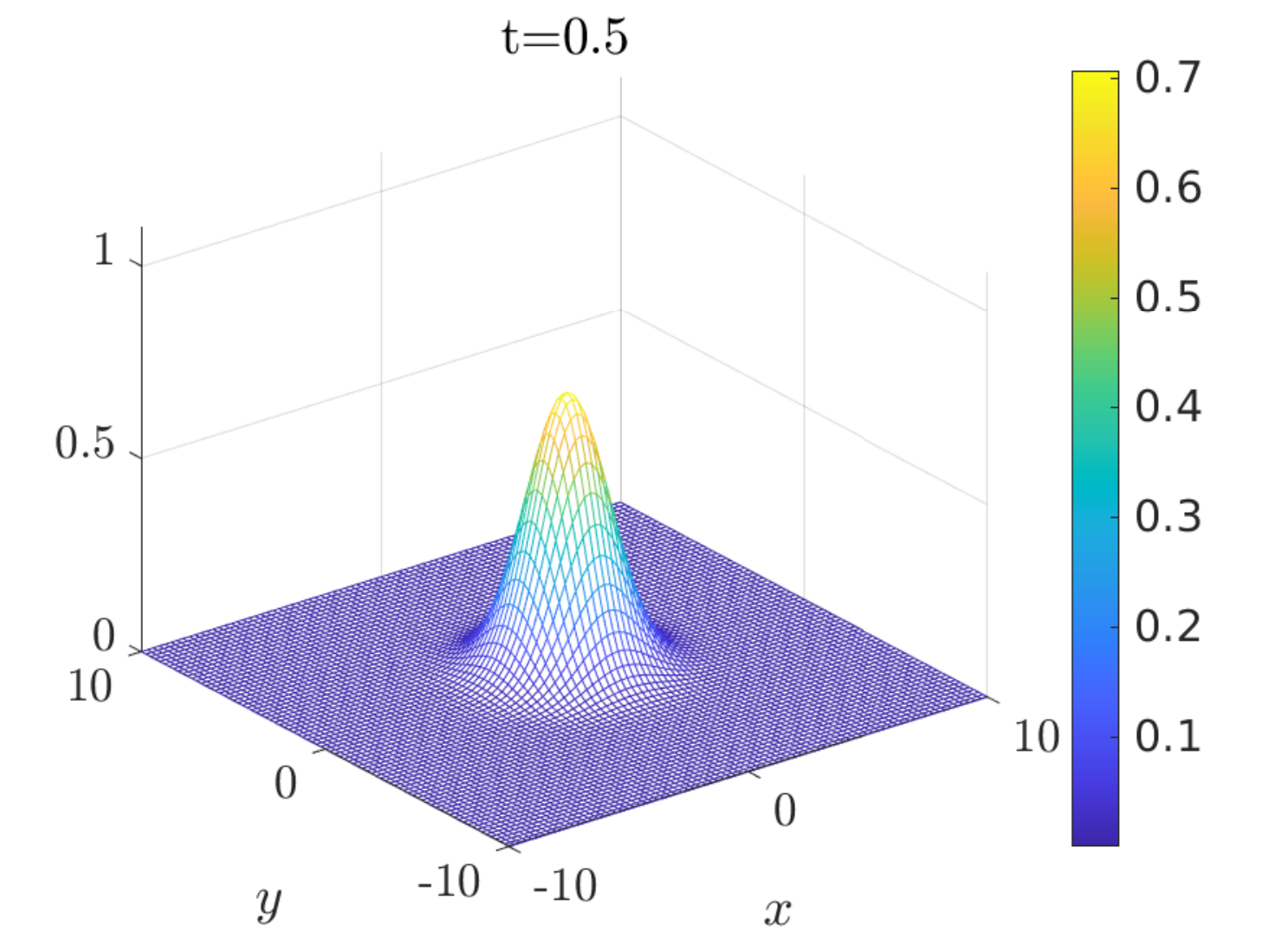}
	\includegraphics[width=4.95cm]{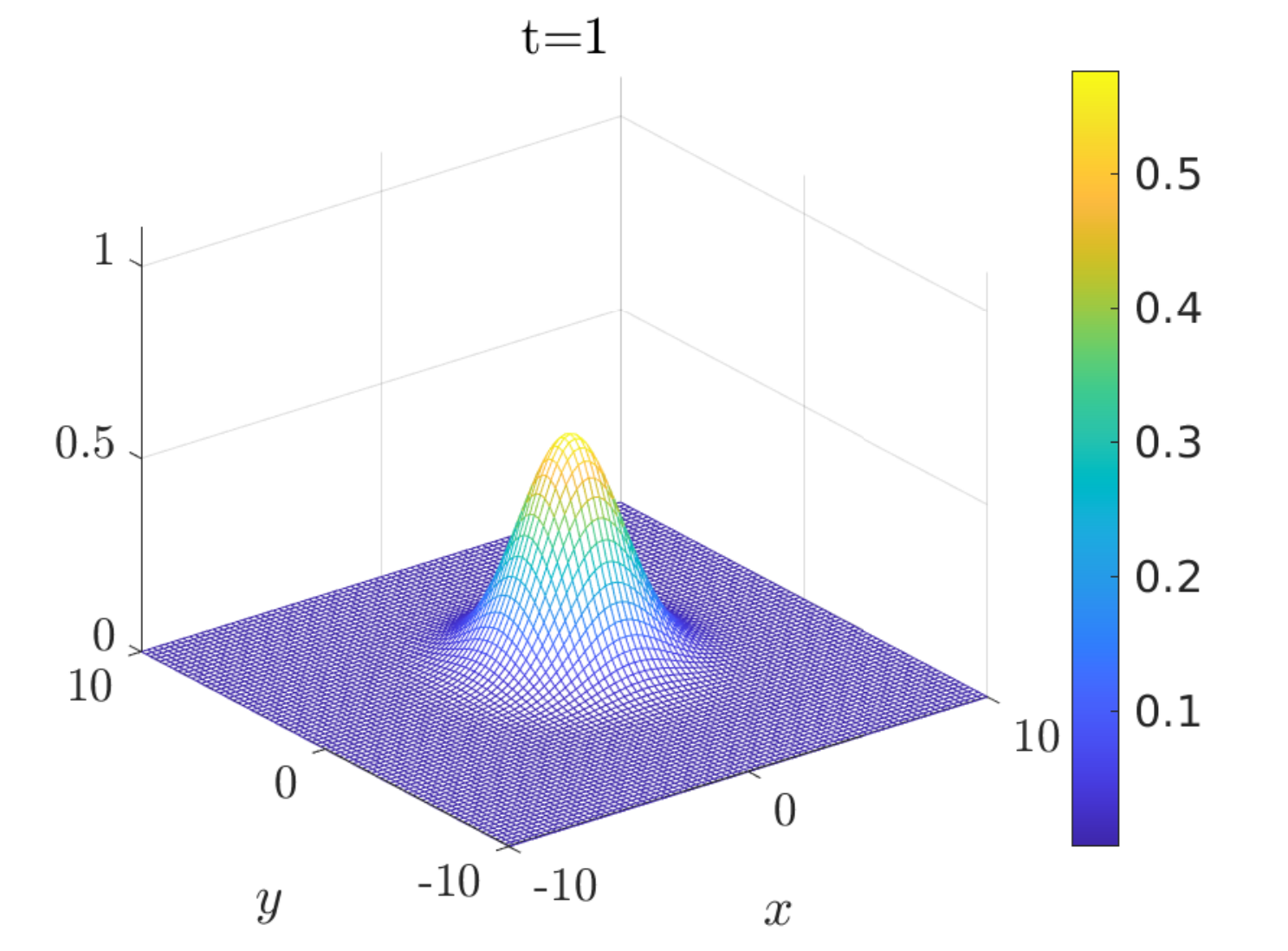}	
	\caption{Test 3: {(\em Non-linear convection diffusion)}. Numerical solution of model \eqref{eq:AD} on the domain $[-10,10]$ up to final time $T=1$. Space step $\Delta x=\Delta y= 0.1$ and $\Delta t =\Delta x/2$. Time integration is performed with LM scheme SSP3-BDF4 and 4th order central difference for convective and diffusion operator.} \label{fig:AD}
\end{figure}

In Figure \ref{fig:AD} we report the numerical solution at different times from $t=0$ to $t = 1$.

\section{Conclusions}
We derived high order semi-implicit schemes based on linear multistep methods. The schemes have been constructed following the approach recently introduced in \cite{BFR} for Runge-Kutta methods. The resulting time discretizations have a predictor-corrector structure and, compared with Runge-Kutta methods, do not require additional order conditions so that they can easily reach high order accuracy. Numerical tests for schemes up to fifth order accuracy have been presented in the case of nonlinear reaction-diffusion and convection-diffusion problems.  

%


%
 \section*{Acknowledgments}

This research is partially supported by PRIN (2019-2021): "Innovative numerical methods for evolutionary partial differential equations and applications" and by INdAM-GNCS (2019) grant: "Approssimazione numerica di problemi di natura iperbolica ed applicazioni".



\end{document}